\documentclass[12pt,draft,reqno]{amsart}

\usepackage[cp1251]{inputenc}
\usepackage[T2A]{fontenc}
\usepackage[english,russian,ukrainian]{babel}
\usepackage{amsmath,amsfonts,amssymb}
\usepackage{geometry}
\author[В. М. Лось, В. А. Михайлець, О. О. Мурач]
{В. М. Лось (V. M. Los), В. А. Михайлець (V. A. Mikhailets), \\ О. О. Мурач (A. A. Murach)}

\title[Параболiчнi задачi у просторах Хермандера]
{Загальнi параболiчнi початково-крайовi задачi \\ у просторах Хермандера \\
\vspace{0.5cm}
General parabolic initial-boundary value problems \\ in H\"ormander spaces}

\address{В. М. Лось, Нацiональний технiчний унiверситет України "КПI",
проспект Перемоги, 37, Kиїв, Україна, 03056}

\email{v$\underline{\ }$los@yahoo.com}

\address{В. А. Михайлець, Iнститут математики НАН України,
вул. Терещенкiвська 3, Kиїв, Україна, 01004}

\email{mikhailets@imath.kiev.ua}

\address{О. О. Мурач, Iнститут математики НАН України,
вул. Терещенкiвська 3, Kиїв, Україна, 01004}

\email{murach@imath.kiev.ua}


\date{20/10/2016}


\begin{document}

\maketitle

\begin{abstract}
We investigate a general nonhomogeneous parabolic initial-boundary value problem
in some anisotropic H\"ormander inner product spaces.
We prove that the operators corresponding to this problem are isomorphisms
between appropriate H\"ormander spaces.
As an application of this result, we establish a theorem on the local increase
in regularity of solutions to the problem.
\end{abstract}

\vspace{1cm}

\textbf{Вступ.} Cучасна теорія загальних параболічних початково--крайових задач
розроблена для класичних шкал функціональних просторів Гельдера--Зігмунда
і Соболєва [1 -- 7]. Центральний результат цієї теорії -- теореми про  коректну
розв'язність (за Адамаром) цих задач у підходящих парах вказаних просторів.
Ці теореми мають важливі застосування до дослідження властивостей регулярності розв'язків
параболічної задачі, властивостей її функції Гріна тощо.

В останній час знаходять важливі застосування у теорії диференціальних рівнянь
із частинними похідними функціональні простори узагальненої гладкості [8 -- 15].
Для них показником регулярності розподілів служить не числовий, а досить загальний
функціональний параметр, залежний від частотних змінних (дуальних до просторових
відносно перетворення Фур'є). Теорія цих просторів була започаткована Л. Хермандером
в монографії \cite{Hermander63}, де було введено нормовані простори
$$\mathcal{B}_{p,\mu}:=\bigl\{ w\in\mathcal{S}'(\mathbb{R}^{k})
\,:\, \mu(\xi)\widehat{w}(\xi)\in L_p(\mathbb{R}^{k},\,d\xi)\bigr\}.
$$
Тут $1\leq p\leq\infty$, $\mu:\mathbb{R}^{k}\rightarrow(0,\infty)$ є вагова
функція, а $\widehat{w}$-- перетворення Фур'є розподілу $w$.

Л. Хермандер дав застосування цих просторів до дослідження регулярності розв'язків
важливих класів лінійних рівнянь у частинних похідних; серед них еліптичні і параболічні
рівняння. Втім довгий час ці простори не знаходили застосування у дослідженні
крайових і початково--крайових задач для цих рівнянь.

Недавно В.~А.~Михайлець і О.~О.~Мурач [14, 15, 17 -- 22]
побудували теорію розв'язності загальних еліптичних крайових
задач у гільбертових шкалах ізотропних просторів Хермандера
$H^{s;\varphi}:=\mathcal{B}_{2,\mu}$, де показник регулярності є функція вигляду
$$
\mu(\xi):=(1+|\xi|^{2})^{s/2}\varphi((1+|\xi|^{2})^{1/2}).
$$
Тут числовий параметр $s$ дійсний, а функціональний параметр $\varphi$ повільно змінний на
нескінченності за Й.~Карамата \cite{Seneta85}.
Ця теорія грунтується на методі інтерполяції з функціональним параметром гільбертових просторів,
який встановлює зв'язок між просторами Соболєва і цими просторами Хермандера.

З огляду на це є перспективним застосування вказаного методу інтерполяції до побудови
теорії розв'язності параболічних початково--крайових задач в анізотропних аналогах
просторів $H^{s;\varphi}$. У наших роботах [24 -- 31] ця теорія розроблена
для мішаних задач з однорідними початковими даними Коші для параболічних рівнянь довільного
порядку та  задач з неоднорідними початковими умовами для рівнянь другого порядку.

Мета цієї роботи -- встановити теореми про коректну розв'язність і локальну регулярність
розв'язків початково-крайової задачі для довільного лінійного $2b$-параболіч\-ного рівняння
з неоднорідними початковими умовами. Першу із згаданих теорем було анонсовано
у \cite{LosMurach14Dop6}.

\textbf{1. Постановка задачі.}
Нехай довільно задані ціле
число $n\geq2$, дійсне число $\tau>0$
і обмежена область $G\subset\mathbb{R}^{n}$ з
нескінченно гладкою межею $\Gamma:=\partial G$.
Позначимо $\Omega:=G\times(0,\tau)$ --- відкритий циліндр в $\mathbb{R}^{n+1}$,
$S:=\Gamma\times(0,\tau)$~--- його бічна поверхня.
Тоді $\overline{\Omega}:=\overline{G}\times[0,\tau]$ і
$\overline{S}:=\Gamma\times[0,\tau]$ є замикання
$\Omega$ і $S$ відповідно.

Розглянемо у циліндрі $\Omega$ таку параболічну початково--крайову задачу
\begin{equation}\label{16f1}
\begin{gathered}
A(x,t,D_x,\partial_t)u(x,t) \equiv\sum_{|\alpha|+2b\beta\leq
2m}a^{\alpha,\beta}(x,t)\,D^\alpha_x\partial^\beta_t
u(x,t)=f(x,t)\\
\mbox{для всіх}\quad x\in G\quad\mbox{і}\quad t\in(0,\tau);
\end{gathered}
\end{equation}
\begin{equation}\label{16f2}
\begin{gathered}
B_{j}(x,t,D_x,\partial_t)u(x,t)\big|_{S} \equiv\sum_{|\alpha|+2b\beta\leq m_j}
b_{j}^{\alpha,\beta}(x,t)\,D^\alpha_x\partial^\beta_t u(x,t)\big|_{S}=g_{j}(x,t)\\
\mbox{для всіх}\quad x\in\Gamma,\quad t\in(0,\tau)\quad\mbox{і}\quad
j\in\{1,\dots,m\};
\end{gathered}
\end{equation}
\begin{equation}\label{16f3}
\partial^k_t u(x,t)
\big|_{t=0}=h_k(x)\quad\mbox{для всіх}\quad x\in G\quad\mbox{і}\quad
k\in\{0,\ldots,\varkappa-1\}.
\end{equation}

Тут $b$, $m$ і всі $m_j$ є довільно задані цілі числа, такі, що $m\geq b\geq1$, $\varkappa:=m/b\in\mathbb{Z}$ і $m_j\geq0$. Число $2b$ називається параболічною вагою цієї задачі.
Усі коефіцієнти лінійних диференціальних виразів $A:=A(x,t,D_x,\partial_t)$ і $B_{j}:=B_{j}(x,t,D_x,\partial_t)$, де $j\in\{1,\dots,m\}$ вважаємо нескінченно гладкими комплекснозначними функціями, заданими на $\overline{\Omega}$ і $\overline{S}$ відповідно; тобто кожна
\begin{equation*}
a^{\alpha,\beta}\in C^{\infty}(\overline{\Omega}):=
\bigl\{w\!\upharpoonright\overline{\Omega}\!:\,w\in C^{\infty}(\mathbb{R}^{n+1})\bigr\}
\end{equation*}
і кожна
\begin{equation*}
b_{j}^{\alpha,\beta}\in C^{\infty}(\overline{S}):=
\bigl\{v\!\upharpoonright\overline{S}\!:\,v\in C^{\infty}(\Gamma\times\mathbb{R})\bigr\}.
\end{equation*}

Використовуємо такі позначення
$D^\alpha_x:=D^{\alpha_1}_{1}\dots D^{\alpha_n}_{n}$, де $D_{k}:=i\,\partial/\partial{x_k}$ і $\partial_t:=\partial/\partial t$
для частинних похідних функцій, що залежать від $x=(x_1,\ldots,x_n)\in\mathbb{R}^{n}$ і $t\in\mathbb{R}$.
Тут $i$ це уявна одиниця, $\alpha=(\alpha_1,...,\alpha_n)$ є мультиіндекс, і $|\alpha|:=\alpha_1+\cdots+\alpha_n$.
У формулах \eqref{16f1} і \eqref{16f2} та їх аналогах підсумовування ведеться за цілими невід'ємними індексами
$\alpha_1,...,\alpha_n$ і $\beta$, які задовольняють умову, вказану під знаком суми.
Як звичайно, $\xi^{\alpha}:=\xi_{1}^{\alpha_{1}}\ldots\xi_{n}^{\alpha_{n}}$ для $\xi:=(\xi_{1},\ldots,\xi_{n})\in\mathbb{C}^{n}$.

Нагадаємо \cite[\S~9, п.~1]{AgranovichVishik64}, що початково--крайова задача
\eqref{16f1}--\eqref{16f3} називається параболічною у  циліндрі $\Omega$, якщо
виконуються такі дві умови.

\emph{Умова $1$.} Для довільних $x\in\overline{G}$, $t\in[0,\tau]$,
$\xi\in\mathbb{R}^{n}$ і $p\in\mathbb{C}$, де $\mathrm{Re}\,p\geq0$, правильно
\begin{equation*}
A^{\circ}(x,t,\xi,p)\equiv\sum_{|\alpha|+2b\beta=2m} a^{\alpha,\beta}(x,t)\,\xi^\alpha
p^{\beta}\neq0\quad\mbox{за умови}\quad|\xi|+|p|\neq0.
\end{equation*}

Для формулювання умови~2 довільно виберемо точку $x\in\Gamma$, дійсне число $t\in[0,\tau]$, дотичний вектор
$\xi\in\mathbb{R}^{n}$ до межи $\Gamma$ у точці $x$ та число $p\in\mathbb{C}$, де
$\mathrm{Re}\,p\geq0$, такі, що $|\xi|+|p|\neq0$. Нехай $\nu(x)$ є ортом
внутрішньої нормалі до межи $\Gamma$ у точці $x$. З умови~1 та нерівності $n\geq2$ випливає, що
многочлен $A^{\circ}(x,t,\xi+\zeta\nu(x),p)$ змінної $\zeta\in\mathbb{C}$ має
рівно $m$ коренів $\zeta^{+}_{j}(x,t,\xi,p)$, $j=\nobreak1,\ldots,m$, с додатною
уявною частиною і $m$ коренів з від'ємною уявною частиною (з урахуванням їх кратності).

\emph{Умова $2$.} При кожному такому виборі $x$, $t$, $\xi$ та $p$ многочлени
$$
B_{j}^{\circ}(x,t,\xi+\zeta\nu(x),p)\equiv\sum_{|\alpha|+2b\beta=m_{j}}
b_{j}^{\alpha,\beta}(x,t)\,(\xi+\zeta\nu(x))^{\alpha}\,p^{\beta},\quad j=1,\dots,m,
$$
змінної $\zeta\in\mathbb{C}$ лінійно незалежні по модулю многочлена
$$
\prod_{j=1}^{m}(\zeta-\zeta^{+}_{j}(x,t,\xi,p)).
$$
Відмітимо, що умова~1 є умовою $2b$-параболічності за І.~Г.~Петровським
\cite{Petrovskii38} диференціаль\-ного рівняння $Au=f$ у замкнутому
циліндрі $\overline{\Omega}$, а умова~2 виражає той факт, що система
крайових диференціальних операторів $\{B_{1},\ldots,B_{m}\}$ накриває
диференціальний оператор $A$ на бічній поверхні $\overline{S}$ цього циліндра.

\textbf{2. Функціональні простори.}
Задачу
\eqref{16f1}--\eqref{16f3} будемо досліджувати у шкалах гільбертових фун\-кці\-ональ\-них
просторів $H^{\mu}:=\mathcal{B}_{2,\mu}$, що були введені Л.~Хермандером у \cite[п.~2.2]{Hermander63}. Згодом ці простори дослідили також та Л.~Р.~Волєвич і Б.~П.~Панеях
\cite[\S~2,~3]{VolevichPaneah65}. Показником регулярності функцій (або розподілів),
що утворюють простір $H^{\mu}(\mathbb{R}^{k})$, де ціле
$k\geq1$, є вимірна за Борелем функція
$\mu:\mathbb{R}^{k}\rightarrow(0,\infty)$, яка задовольняє таку умову:
існують додатні числа $c$ та $l$ такі, що
$$
\frac{\mu(\xi)}{\mu(\eta)}\leq
c\,(1+|\xi-\eta|)^{l}\quad\mbox{для довільних}\quad \xi,\eta\in\mathbb{R}^{k}.
$$

За означенням, комплексний лінійний простір $H^{\mu}(\mathbb{R}^{k})$ складається з усіх
повільно зростаючих розподілів $w\in\mathcal{S}'(\mathbb{R}^{k})$, перетворення
Фур'є $\widehat{w}$ яких є локально інтегровними за Лебегом функціями, що
задовольняють умову
\begin{equation*}
\int\limits_{\mathbb{R}^{k}}\mu^{2}(\xi)\,|\widehat{w}(\xi)|^{2}\,d\xi
<\infty.
\end{equation*}
(У роботі усі функції та розподіли вважаються комплекснозначними.)
У просторі $H^{\mu}(\mathbb{R}^{k})$ означений скалярний добуток
за формулою
\begin{equation*}
(w_1,w_2)_{H^{\mu}(\mathbb{R}^{k})}=
\int\limits_{\mathbb{R}^{k}}\mu^{2}(\xi)\,\widehat{w_1}(\xi)\,
\overline{\widehat{w_2}(\xi)}\,d\xi,
\end{equation*}
де $w_1,w_2\in H^{\mu}(\mathbb{R}^{k})$. Цей скалярний добуток породжує
норму
$$
\|w\|_{H^{\mu}(\mathbb{R}^{k})}:=(w,w)^{1/2}_
{H^{\mu}(\mathbb{R}^{k})}.
$$
Простір $H^{\mu}(\mathbb{R}^{k})$ є гільбертовим і сепарабельним відносно введеного у ньому скалярного добутку. Цей простір є неперервно вкладеним у $\mathcal{S}'(\mathbb{R}^{k})$,
а множина $C^{\infty}_{0}(\mathbb{R}^{k})$ є щільною в ньому
\cite[п. 2.2]{Hermander63}.

Нам знадобиться версія простору $H^{\mu}(\mathbb{R}^{k})$ для
довільної відкритої множини $V\neq\varnothing$. Лінійний простір
$H^{\mu}(V)$ складається, за означеннням, із звужень $u=w\!\upharpoonright\!V$ всіх
розподілів $w\in H^{\mu}(\mathbb{R}^{k})$ на множину $V$. У цьому просторі
задана норма за формулою
\begin{equation}\label{16f68}
\|u\|_{H^{\mu}(V)}:= \inf\bigl\{\|w\|_{H^{\mu}(\mathbb{R}^{k})}:\,w\in
H^{\mu}(\mathbb{R}^{k}),\;u=w\!\upharpoonright\!V\bigr\}.
\end{equation}
Іншими словами,
$H^{\mu}(V)$ є фактор-простором простору
$H^{\mu}(\mathbb{R}^{k})$ за його підпростором
\begin{equation}\label{16f69}
H^{\mu}(\mathbb{R}^{k},V):=\bigl\{w\in
H^{\mu}(\mathbb{R}^{k}):\,
w=0\;\,\mbox{на}\;\,V\bigr\}.
\end{equation}
Тому простір $H^{\mu}(V)$ є гільбертовим і сепарабельним. Норма
\eqref{16f68} породжена скалярним добутком
$$
(u_{1},u_{2})_{H^{\mu}(V)}:= (w_{1}-\Upsilon
w_{1},w_{2}-\Upsilon w_{2})_{H^{\mu}(\mathbb{R}^{k})},
$$
де $w_{j}\in H^{\mu}(\mathbb{R}^{k})$, $w_{j}=u_{j}$ у $V$
для кожного $j\in\{1,\,2\}$. Тут $\Upsilon$ є ортогональним проектором простору
$H^{\mu}(\mathbb{R}^{k})$ на його підпростір \eqref{16f69}.
У просторі $H^{\mu}(V)$ щільна множина
$$
C^{\infty}_{0}(\overline{V}):=\bigl\{w\upharpoonright_{\overline{V}}
\,\,:\,w\in C^{\infty}_{0}(\mathbb{R}^{k})\bigr\}.
$$

Нехай задане число $\gamma>0$. Надалі будемо використовувати
показники регулярності вигляду
\begin{equation}\label{16f70}
\mu_{s,\varphi}(\xi',\xi_{k}):=
\mu(\xi',\xi_{k}):=\bigl(1+|\xi'|^2+|\xi_{k}|^{2\gamma}\bigr)^{s/2}
\varphi\bigl((1+|\xi'|^2+|\xi_{k}|^{2\gamma})^{1/2}\bigr),
\end{equation}
де $\xi'\in\mathbb{R}^{k-1}$ та $\xi_{k}\in\mathbb{R}$ є аргументами функції $\mu$.
Тут числовий параметр $s$ є дійсним, а функціональний параметр $\varphi$
пробігає клас $\mathcal{M}$.

За означенням, клас $\mathcal{M}$ складається з усіх вимірних за Борелем функцій $\varphi:[1,\infty)\rightarrow(0,\infty)$, які задовольняють такі дві умови:

а) обидві функції $\varphi$ та $1/\varphi$ обмежені на кожному відрізку $[1,c]$, де
$1<c<\infty$;

б) функція $\varphi$ повільно змінюється за Й.~Карамата на нескінченності, а саме,
$\varphi(\lambda r)/\varphi(r)\rightarrow 1$ при $r\rightarrow\infty$ для кожного
$\lambda>0$.

Теорія повільно змінних функцій (на нескінченності) викладена, наприклад, у
монографії \cite{Seneta85}. Їх важливим прикладом є  функції вигляду
\begin{equation*}
\varphi(r):=(\log r)^{\theta_{1}}\,(\log\log r)^{\theta_{2}} \ldots
(\,\underbrace{\log\ldots\log}_{k\;\mbox{\small{разів}}}r\,)^{\theta_{k}}
\quad\mbox{при}\quad r\gg1,
\end{equation*}
де параметри $k\in\mathbb{N}$ та
$\theta_{1},\theta_{2},\ldots,\theta_{k}\in\mathbb{R}$ є довільними.

Відмітимо, що простори, які введемо нижче будуть потрібні нам
лише у випадку $\gamma=1/(2b)$. Але природньо їх ввести для
довільного $\gamma>0$.

Нехай $s\in\mathbb{R}$ і $\varphi\in\mathcal{M}$.
Розв'язки $u$ початково--крайової задачі \eqref{16f1}--\eqref{16f3}
та праві частини $f$ рівняння \eqref{16f1} будемо
розглядати у анізотропних гільбертових функціональних просторах Хермандера
$H^{s,s\gamma;\varphi}(\Omega):=H^{\mu}(\Omega)$ де показник $\mu$
визначений формулою~\eqref{16f70}, у якій $k:=n+1$.

Якщо $\varphi(r)\equiv1$, то $H^{s,s\gamma;\varphi}(\Omega)$ стає анізотропним
гільбертовим простором Соболева порядку $(s,s\gamma)$; позначимо його через
$H^{s,s\gamma}(\Omega)$. Тут $s$~--- показник регулярності розподілу
$u=u(x,t)$ по просторовій змінній $x\in\Omega$, а $s\gamma$~--- показник
регулярності по часовій змінній $t\in(0,\tau)$. В загальному випадку, коли
$\varphi\in\mathcal{M}$ є довільною, правильні неперервні і щільні вкладення
\begin{equation}\label{16f71}
H^{s_{1},s_{1}\gamma}(\Omega)\hookrightarrow
H^{s,s\gamma;\varphi}(\Omega)\hookrightarrow
H^{s_{0},s_{0}\gamma}(\Omega)\quad\mbox{при}\quad s_{0}<s<s_{1}.
\end{equation}
Справді, оскільки функція $\varphi\in\mathcal{M}$, то існують додатні числа
$c_0$ і $c_1$ такі, що
$$
c_0\,r^{s_0-s}\leq\varphi(r)\leq c_1\,r^{s_1-s}\quad\mbox{для
довільного}\quad r\geq1
$$
(див., наприклад, \cite[п.1.5,$1^0$]{Seneta85}). Тоді
\begin{equation*}
c_0\,\mu_{s_0,1}(\xi',\xi_{n+1})\leq\mu_{s,\varphi}(\xi',\xi_{n+1})\leq
c_1\,\mu_{s_1,1}(\xi',\xi_{n+1})
\end{equation*}
для довільних $\xi'\in\mathbb{R}^{n}$ та $\xi_{n+1}\in\mathbb{R}$.
Звідси зразу випливають неперервні вкладення просторів \eqref{16f71}.
Ці вкладення щільні, оскільки множина $C^{\infty}_{0}(\overline{\Omega})$
щільна в усіх цих просторах.

Розглянемо клас гільбертових функціональних просторів
\begin{equation}\label{16f72}
\bigl\{H^{s,s\gamma;\varphi}(\Omega):\,
s\in\mathbb{R},\,\varphi\in\mathcal{M}\,\bigr\}.
\end{equation}
Вкладання \eqref{16f71} показують, що у \eqref{16f72} функціональний параметр $\varphi$
визначає додаткову гладкість по відношенню до основної анізотропної
$(s,s\gamma)$-гладкості. Якщо $\varphi(r)\rightarrow\infty$ (або
$\varphi(r)\rightarrow0$) при $r\rightarrow\infty$, то $\varphi$ визначає
позитивну (або негативну) додаткову гладкість. Інакше кажучи, $\varphi$
уточнює основну гладкість $(s,s\gamma)$.
Тут $\gamma>0$ виконує роль параметра анізотропії
просторів, що утворюють цю шкалу.

Нам знадобляться також анізотропні простори Хермандера, задані на
бічній поверхні $S=\Gamma\times(0,\tau)$ циліндра $\Omega$.
До них будуть належати праві частини $g_j$ крайових
умов \eqref{16f2}.  Означимо ці простори використовуючи
спеціальні локальні карти на $S$ (див. \cite[п.1]{Los15NK2}).

Нехай $s>0$ і $\varphi\in\mathcal{M}$.
Попередньо для відкритої смуги $\Pi:=\mathbb{R}^{n-1}\times(0,\tau)$ розглянемо гільбертові
простори $H^{s,s\gamma;\varphi}(\Pi):=H^{\mu}(\Pi)$, де показник $\mu$ визначений формулою
\eqref{16f70}, у якій $k:=n$.
Довільно виберемо скінченний атлас із
$C^{\infty}$-структури на замкненому многовиді~$\Gamma$. Нехай цей атлас утворений локальними
картами $\nobreak{\theta_{j}:\mathbb{R}^{n-1}\leftrightarrow \Gamma_{j}}$, де
$j=1,\ldots,\lambda$. Тут відкрити множини $\Gamma_{1},\ldots,\Gamma_{\lambda}$
складають покриття многовиду $\Gamma$. Окрім цього, довільно виберемо функції
$\chi_{j}\in C^{\infty}(\Gamma)$, $j=1,\ldots,\lambda$, такі, що
$\mathrm{supp}\,\chi_{j}\subset\Gamma_{j}$ і $\sum_{j=1}^{\lambda}\chi_{j}\equiv1$
на $\Gamma$.

За означенням,
лінійний простір $H^{s,s\gamma;\varphi}(S)$ складається з усіх
функцій $v\in L_2(S)$ на многовиді $S$ таких, що для кожного
номеру $j\in\{1,\ldots,\lambda\}$ функція
$$
v_{j}(x,t):=\chi_{j}(\theta_{j}(x))\,v(\theta_{j}(x),t)
$$
аргументів $x\in\mathbb{R}^{n-1}$ і $t\in(0,\tau)$ належить до $H^{s,s\gamma;\varphi}(\Pi)$.
У просторі $H^{s,s\gamma;\varphi}(S)$ означений скалярний добуток за формулою
\begin{equation*}
(v,g)_{H^{s,s\gamma;\varphi}(S)}:=
\sum_{j=1}^{\lambda}\,
(v_{j},g_{j})_{H^{s,s\gamma;\varphi}(\Pi)},
\end{equation*}
де $v,g\in H^{s,s\gamma;\varphi}(S)$.
Він породжує норму
$$
\|v\|_{H^{s,s\gamma;\varphi}(S)}:=(v,v)^{1/2}_{H^{s,s\gamma;\varphi}(S)}.
$$

Відмітимо, що простір $H^{s,s\gamma;\varphi}(S)$ означений за допомогою
спеціальних локальних карт на $S$
\begin{equation*}
\theta_{j}^*: \Pi=\mathbb{R}^{n-1}\times(0,\tau)\leftrightarrow\Gamma_{j}\times(0,\tau)
\quad\mbox{для кожного}\quad j\in\{1,\ldots,\lambda\},
\end{equation*}
де покладаємо $\theta_{j}^*(x,t):=(\theta_{j}(x),t)$ для усіх
$x\in\mathbb{R}^{n-1}$,\, $t\in(0,\tau)$.
Цей простір є повним (гільбертовим) та не залежить з точністю до еквівалентності норм
від вибору локальних карт і розбиття одиниці на $\Gamma$ \cite[теорема 1]{Los15NK2}.

Нам також знадобляться ізотропні простори Хермандера $H^{s;\varphi}(V)$, де
$s\in\mathbb{R}$, $\varphi\in\mathcal{M}$, а $V$~--- довільна відкрита непорожня множина.
За означенням,
$H^{s;\varphi}(V)$ є гільбертів простір
$H^{\mu}(V)$, де
\begin{equation*}
\mu(\xi)=\bigl(1+|\xi|^2\bigr)^{s/2}
\varphi\bigl((1+|\xi|^2)^{1/2}\bigr).
\end{equation*}
Тут $\xi\in\mathbb{R}^{k}$ є аргументом функції $\mu$.
Оскільки функція $\mu$ є радіальною (залежить лише від
$|\xi|$), то простір $H^{s;\varphi}(V)$ є ізотропним.
Ми будемо використовувати простори $H^{s;\varphi}(V)$ у
випадках $V=\mathbb{R}^{k}$ і $V=G$.
До просторів $H^{s;\varphi}(G)$ будуть належати
праві частини $h_k$ початкових умов \eqref{16f3}.

Окрім того, як допоміжні, нам будуть потрібні простори $H^{s;\varphi}(\Gamma)$.
Подібно до просторів на $S$, означимо їх за допомогою локальних карт
$\theta_{j}$, $j=1,\ldots,\lambda$, вказаних вище.
Нехай $s\in\mathbb{R}$ і $\varphi\in\mathcal{M}$.
За означенням,
лінійний простір $H^{s;\varphi}(\Gamma)$ складається з усіх
розподілів $w\in\mathcal{D}'(\Gamma)$ на многовиді $\Gamma$ таких, що для кожного
номеру $j\in\{1,\ldots,\lambda\}$ розподіл
$w_{j}(x):=\chi_{j}(\theta_{j}(x))\,w(\theta_{j}(x))$ аргумента
$x\in\mathbb{R}^{n-1}$ належить до $H^{s;\varphi}(\mathbb{R}^{n-1})$.
У просторі $H^{s;\varphi}(\Gamma)$ означений скалярний добуток за формулою
\begin{equation*}
(w,g)_{H^{s;\varphi}(\Gamma)}:=
\sum_{j=1}^{\lambda}\,
(w_{j},g_{j})_{H^{s;\varphi}(\mathbb{R}^{n-1})},
\end{equation*}
де $w,g\in H^{s;\varphi}(\Gamma)$.
Він породжує норму
$$
\|w\|_{H^{s;\varphi}(\Gamma)}:=(w,w)^{1/2}_{H^{s;\varphi}(\Gamma)}.
$$
Простір $H^{s;\varphi}(\Gamma)$ є гільбертовим та не залежить з точністю до еквівалентності норм
від вибору локальних карт і розбиття одиниці на $\Gamma$ \cite[теорема 2.3]{MikhailetsMurach14}.

Ізотропні простори $H^{s;\varphi}$ виділили і систематично використовували
В.~А.~Михайлець та О.~О.~Мурач у теорії еліптичних крайових
задач~\cite{12BJMA2, MikhailetsMurach14}.

Якщо $\varphi\equiv1$, то означені вище простори стають соболевськими
просторами (анізотропними або ізотропними).
У цьому випадку будемо опускати індекс $\varphi$ у позначеннях
цих і введених нижче просторів.

\textbf{3. Основні результати}
становлять теорему про ізоморфізми для параболічної задачі \eqref{16f1}--\eqref{16f3}
у введених вище просторах Хермандера і застосування цієї теореми до дослідження
локальної регулярності узагальнених розв'язків задачі. Сформулюємо їх.

Нехай $\sigma_0$ є найменше ціле число, таке, що
$$
\sigma_0\geq2m,\quad\sigma_0\geq m_j+1\;\;\mbox{для всіх}\;\;j\in\{1,\ldots,m\}.
$$
Відмітимо, якщо $m_j\leq2m-1$ для всіх $j\in\{1,\ldots,m\}$, тоді $\sigma_0=2m$.

Попередньо введемо функціональний простір,
до якого буде належати вектор--функція $(f,g_1,\dots,g_m,h_0,\dots,h_{\varkappa-1})$ правих частин задачі.
Елементи цього простору задовольняють деякі умови узгодження
(див. \cite[\S11]{AgranovichVishik64} або \cite[гл.4,\S5]{LadyzhenskajaSolonnikovUraltzeva67}).
Ці умови полягають в тому, що похідні $\partial^k_t u(x,t)\big|_{t=0}$, які можна обчислити з
параболічного рівняння \eqref{16f1} та початкових умов \eqref{16f3}, повинні
задовольняти при $x\in\Gamma$ крайові умови \eqref{16f2} та співвідношення, що утворюються
в результаті диференціювання цих крайових умов по змінній $t$.
Розглянемо ці умови, спираючись на соболєвську теорію параболічних задач ($\varphi\equiv1$).

Пов'яжемо із задачею \eqref{16f1}, \eqref{16f2}, \eqref{16f3} лінійне відображення
\begin{equation}\label{16f4}
u\mapsto\Lambda
u:=\bigl(Au,B_1u,\ldots,B_mu,u\!\upharpoonright\!\overline{G},\ldots,
(\partial^{\varkappa-1}_tu)\!\upharpoonright\!\overline{G}\bigr),\quad u\in
C^{\infty}(\overline{\Omega}).
\end{equation}
Нехай $s\geq\sigma_0$. Відображення \eqref{16f4} однозначно продовжується (за неперервністю) до
обмеженого оператора
\begin{equation*}
\begin{gathered}
\Lambda:\,H^{s,s/(2b)}(\Omega)\; \rightarrow\;
\mathcal{H}^{s-2m,(s-2m)/(2b)}:=
H^{s-2m,(s-2m)/(2b)}(\Omega)\oplus\\
\oplus\bigoplus_{j=1}^{m}H^{s-m_j-1/2,(s-m_j-1/2)/(2b)}(S)
\oplus\bigoplus_{k=0}^{\varkappa-1}H^{s-2bk-b}(G).
\end{gathered}
\end{equation*}
Це означає, що для кожної функції $u=u(x,t)$ з соболєвського простору
$H^{s,s/(2b)}(\Omega)$ означені праві частини задачі
\begin{equation*}
\begin{gathered}
f\in H^{s-2m,(s-2m)/(2b)}(\Omega),\quad g_j\in H^{s-m_j-1/2,(s-m_j-1/2)/(2b)}(S)\quad
\mbox{і}\quad h_k\in H^{s-2bk-b}(G)\\
\mbox{для всіх}\quad j\in\{1,\dots,m\}\quad\mbox{і}\quad
k\in\{0,\dots,\varkappa-1\}.
\end{gathered}
\end{equation*}

Умови узгодження на функції $f$, $g_j$ і $h_k$
природно виникають так. Для функції $u$ означені сліди
$\partial^{\,k}_t u(\cdot,0)\in H^{s-2bk-b}(G)$ для всіх цілих $k$,
таких, що $0\leq k<s/(2b)-1/2$, і лише таких цілих $k$.
Ці сліди виражаються з рівняння \eqref{16f1} і початкових умов
\eqref{16f3} через функції $f$ і $h_k$ наступним чином.

Умова параболічності у випадку $\xi=0$ і $p=1$ означає, що коефіцієнт $a^{(0,\ldots,0),\varkappa}(x,t)\neq0$
для всіх $x\in\overline{G}$ і $t\in[0,\tau]$. Тому параболічне рівняння \eqref{16f1} можна розв'язати відносно  $\partial^\varkappa_t u(x,t)$; а саме, запишемо
\begin{equation}\label{16f5}
\partial^\varkappa_t u(x,t)=
\sum_{\substack{|\alpha|+2b\beta\leq 2m,\\ \beta\leq\varkappa-1}}
a_{0}^{\alpha,\beta}(x,t)\,D^\alpha_x\partial^\beta_t
u(x,t)+f(x,t)
\end{equation}
для деяких функцій $a_{0}^{\alpha,\beta}\in C^{\infty}(\overline{\Omega})$.
З умов \eqref{16f3}, рівності \eqref{16f5} та рівностей, одержаних з неї шляхом диференціювання
по змінній $t$ необхідну кількість разів отримуємо рекурентні формули для слідів
$\partial^{\,k}_t u(x,0)$:
\begin{equation}\label{16f6}
\begin{split}
&\partial^{\,k}_t u(x,0)=h_k(x)\quad\mbox{для всіх}
\quad k\in\{0,\dots,\varkappa-1\},\\
&\partial^{\,k}_t u(x,0)=
\sum_{\substack{|\alpha|+2b\beta\leq 2m,\\ \beta\leq\varkappa-1}}
\sum\limits_{q=0}^{k-\varkappa}
\binom{k-\varkappa}{q}\partial^{\,k-\varkappa-q}_t
a_{0}^{\alpha,\beta}(x,0)\,D^\alpha_x\partial^{\beta+q}_t
u(x,0)+\\
&+\partial^{\,k-\varkappa}_t f(x,0)
\quad\mbox{для}\quad k\geq\varkappa;
\end{split}
\end{equation}
тут $x\in G$ -- довільне, а цілі $k$, такі, що $0\leq k<s/(2b)-1/2$.

Окрім того, для кожного $j\in\{1,\dots,m\}$
сліди $\partial^{\,k}_t g_j(\cdot,0)\in H^{s-m_j-1/2-2bk-b}(\Gamma)$ означені
для усіх цілих чисел $k$, таких, що $0\leq k<(s-m_j-1/2-b)/(2b)$, і лише для цих цілих $k$.
Тому, згідно з крайовими умовами \eqref{16f2}, для кожного $j\in\{1,\dots,m\}$ виконуються рівності
\begin{equation}\label{16f7}
(\partial^{\,k}_t B_j u)(\cdot,0)=(\partial^{\,k}_t g_j)(\cdot,0)
\end{equation}
на $\Gamma$ для усіх цих цілих $k$.

Для кожного $j\in\{1,\dots,m\}$ і кожного цілого $k\geq0$ позначимо
\begin{align*}
B_{j,k}(u,\partial_t u,\dots,&\partial^{\,[m_j/(2b)]+k}_t u):=\partial^{\,k}_t B_ju|_{t=0}=\\ \\
&\sum_{|\alpha|+2b\beta\leq m_j}
\sum_{q=0}^{k}\binom{k}{q}
\partial^{\,k-q}_t b^{\alpha,\beta}_j(x,0)
D^\alpha_x\partial^{\,\beta+q}_t u(x,0).
\end{align*}

Тепер підставивши \eqref{16f6} у \eqref{16f7} дістаємо
умови узгодження
\begin{equation}\label{16f8}
\begin{gathered}
\partial^{\,k}_t g_j\!\upharpoonright\!\Gamma=
B_{j,k}(v_0, v_1,\dots,v_{[m_j/(2b)]+k})\!\upharpoonright\!\Gamma \\ \\
\mbox{для всіх}\;\;\;k\in\mathbb{Z}\;\;\;\mbox{таких, що}\;\;\;
0\leq k<\frac{s-m_j-1/2-b}{2b}\;\;\;\mbox{і всіх}\;\;\;
j\in\{1,\dots,m\}.
\end{gathered}
\end{equation}
Тут
\begin{equation*}
B_{j,k}(v_0, v_1,\dots,v_{[m_j/(2b)]+k})=
\sum_{|\alpha|+2b\beta\leq m_j}
\sum_{q=0}^{k}\binom{k}{q}
\partial^{\,k-q}_t b^{\alpha,\beta}_j(x,0)
D^\alpha_x v_{\beta+q}(x),
\end{equation*}
функції $v_l=v_l(\cdot,f,h_0,\dots,h_{\varkappa-1})$ означені за рекурентними формулами
\begin{equation}\label{16f9}
\begin{split}
&v_l(x)=h_l(x)\quad\mbox{для всіх}
\quad l\in\{0,\dots,\varkappa-1\},\\
&v_l(x)=
\sum_{\substack{|\alpha|+2b\beta\leq 2m,\\ \beta\leq\varkappa-1}}
\sum\limits_{q=0}^{l-\varkappa}
\binom{l-\varkappa}{q}\partial^{\,l-\varkappa-q}_t
a_{0}^{\alpha,\beta}(x,0)\,D^\alpha_x v_{\beta+q}(x)+\\
&+\partial^{\,l-\varkappa}_t f(x,0)
\quad\mbox{для}\quad l\geq\varkappa;
\end{split}
\end{equation}
де $x\in G$ -- довільне, а цілі $l$, такі, що $0\leq l<(s-m_j-1/2-b)/(2b)$.

Для довільних параметрів $s\geq\sigma_0$ і $\varphi\in\mathcal{M}$ введемо гільбертів простір $\mathcal{Q}^{s-2m,(s-2m)/(2b);\varphi}$ правих частин задачі
\eqref{16f1}--\eqref{16f3}.
Позначимо
\begin{equation*}
\begin{split}
&\mathcal{H}^{s-2m,(s-2m)/(2b);\varphi}:=\\
&H^{s-2m,(s-2m)/(2b);\varphi}(\Omega)
\oplus\bigoplus_{j=1}^{m}H^{s-m_j-1/2,(s-m_j-1/2)/(2b);\varphi}(S)
\oplus\bigoplus_{k=0}^{\varkappa-1}H^{s-2bk-b;\varphi}(G).
\end{split}
\end{equation*}
Якщо $s\notin\{\sigma_0+r-1/2\,:\,r\in\mathbb{N}\}$, то, за означенням, лінійний простір $\mathcal{Q}^{s-2m,(s-2m)/(2b);\varphi}$ складається з усіх вектор--функцій
$$F=\bigl(f,g_1,\dots,g_m,h_0,\dots,h_{\varkappa-1}\bigr)\in
\mathcal{H}^{s-2m,(s-2m)/(2b);\varphi},$$
які задовольняють умови узгодження \eqref{16f8}. Цей лінійний простір наділяється скалярним добутком і нормою з гільбертового простору
$\mathcal{H}^{s-2m,(s-2m)/(2b);\varphi}$. Простір $\mathcal{Q}^{s-2m,(s-2m)/(2b);\varphi}$ є
повним, тобто гільбертовим.

Цей факт відомий у соболєвському випадку $\varphi\equiv1$  і є наслідком обмеженості у відповідних анізотропних соболєвських просторах диференціальних і крайових операторів, що фігурують у формулах \eqref{16f8} і~\eqref{16f9}. У загальній ситуації, коли $\varphi\in\mathcal{M}$, повнота простору $\mathcal{Q}^{s-2m,(s-2m)/(2b);\varphi}$
є наслідком рівності
$$
\mathcal{Q}^{s-2m,(s-2m)/(2b);\varphi}=
\mathcal{H}^{s-2m,(s-2m)/(2b);\varphi}\cap
\mathcal{Q}^{s-2m-\varepsilon,(s-2m-\varepsilon)/(2b)}.
$$
Тут число $\varepsilon>0$ настільки мале, що
елементи просторів
$$
\mathcal{Q}^{s-2m,(s-2m)/(2b);\varphi}\quad\mbox{і}\quad
\mathcal{Q}^{s-2m-\varepsilon,(s-2m-\varepsilon)/(2b)}
$$
задовольняють одні і ті самі умови узгодження \eqref{16f8}. Права частина цієї рівності є гільбертовим простором відносно скалярного добутку, що є сумою скалярних добутків у гільбертових просторах~--- компонентах перетину. Окрім того, норми у просторах, з'єднаних знаком рівності, еквівалентні на підставі неперервного вкладення
$$
\mathcal{H}^{s-2m,(s-2m)/(2b);\varphi}\hookrightarrow
\mathcal{H}^{s-2m-\varepsilon,(s-2m-\varepsilon)/(2b)},
$$
яке є наслідком формули \eqref{16f71}. Тому є повним і простір у лівій частині цієї рівності. З~останнього вкладення випливає також, що умови узгодження \eqref{16f8} коректно поставлені для довільного вектора $F\in\mathcal{H}^{s-2m,(s-2m)/(2b);\varphi}$, оскільки вони є коректними, якщо цей вектор лежить у ширшому анізотропному соболєвському просторі $\mathcal{H}^{s-2m-\varepsilon,(s-2m-\varepsilon)/(2b)}$.

Якщо $s\in\{\sigma_0+r-1/2\,:\,r\in\mathbb{N}\}$, то
означаємо гільбертів простір  $\mathcal{Q}^{s-2m,(s-2m)/(2b);\varphi}$
за допомогою інтерполяції
\begin{equation}\label{16f10}
\mathcal{Q}^{s-2m,(s-2m)/(2b);\varphi}:=
\bigl[\mathcal{Q}^{s-2m-\varepsilon,(s-2m-\varepsilon)/(2b);\varphi},
\mathcal{Q}^{s-2m+\varepsilon,(s-2m+\varepsilon)/(2b);\varphi}\bigr]_{1/2}.
\end{equation}
Тут число $\varepsilon\in(0,1/2)$ вибрано довільно, а права частина рівності є результатом інтерполяції зазначеної пари гільбертових просторів з числовим параметром $1/2$ (див. п.~4). Означений у такий спосіб простір не залежить з точністю до еквівалентності норм від вибору числа $\varepsilon$
(див. зауваження у п.~6 після доведення теореми~1).

\textbf{Теорема 1.} \it Нехай довільно задані параметри: числовий
$s>\sigma_0$ і функціональний $\varphi\in\nobreak\mathcal{M}$.
Тоді відображення \eqref{16f4} однозначно продовжується (за неперервністю) до ізоморфізму
\begin{equation}\label{16f11}
\Lambda :\,H^{s,s/(2b);\varphi}(\Omega)\leftrightarrow
\mathcal{Q}^{s-2m,(s-2m)/(2b);\varphi}.
\end{equation}
\rm
У соболєвському випадку $\varphi\equiv1$ ця теорема є відомою.
Вона доведена М.~С.~Аграновічем та М.~І.~Вішиком \cite[теорема 12.1]{AgranovichVishik64} у
припущенні $s/(2b)\in\mathbb{Z}$. Це припущення можна зняти, що випливає з результату
М.~В.~Житарашу \cite[теорема 9.1]{Zhitarashu85}.
Його результат охоплює граничний випадок $s=\sigma_{0}$.
У загальному випадку виведемо теорему~1 з цього результату шляхом інтерполяції
з функціональним параметром пар гільбертових просторів Соболєва.
Це буде зроблено у п.6 після того, як дослідимо необхідні інтерполяційні
властивості просторів Хермандера, що з'являються у \eqref{16f11}.

Зі згаданої теореми М.~В.~ Житарашу випливає, що для кожної вектор--функції
\begin{equation}\label{16f12}
(f,g_{1},\dots,g_{m},h_0,\dots,h_{\varkappa-1})\in
\mathcal{Q}^{\sigma_0-2m,(\sigma_0-2m)/(2b)}
\end{equation}
задача \eqref{16f1}--\eqref{16f3} має єдиний розв'язок
$u\in H^{\sigma_0,\sigma_0/(2b)}(\Omega)$.
Таку функцію $u$ називаємо узагальненим розв'язком цієї задачі
із правою частиною \eqref{16f12}.

Обговоримо властивості регулярності цього розв'язку у просторах Хермандера.
Наступний результат негайно випливає з теореми~1.

\noindent\textbf{Наслідок 1.} \it Припустимо, що функція $u\in
H^{\sigma_0,\sigma_0/(2b)}(\Omega)$ є узагальненим розв'язком параболічної задачі
\eqref{16f1}--\eqref{16f3}, праві частини якої задовольняють умову
$$
(f,g_{1},\dots,g_{m},h_0,\dots,h_{\varkappa-1})\in
\mathcal{Q}^{\sigma-2m,(\sigma-2m)/(2b);\varphi}
$$
для деяких дійсного числа $\sigma>\sigma_0$ і функції $\varphi\in\mathcal{M}$.
Тоді $u\in H^{\sigma,\sigma/(2b);\varphi}(\Omega)$.
\rm

Тепер сформулюємо локальний аналог цього результату.
Нехай $U$ є відкритою множиною в $\mathbb{R}^{n+1}$, такою, що $U\cap\Gamma=\varnothing$.
Нехай $\omega:=U\cap\Omega\neq\varnothing$,
$\pi_{1}:=U\cap\partial\Omega$,  $\pi_{2}:=U\cap S$ і $\pi_{3}:=U\cap G$.
Введемо необхідні локальні аналоги просторів
$H^{s,s\gamma;\varphi}(\Omega)$,  $H^{s,s\gamma;\varphi}(S)$ і $H^{s;\varphi}(G)$
з $s>0$, $\gamma=1/(2b)$, і $\varphi\in\mathcal{M}$.

Позначимо через $H^{s,s\gamma;\varphi}_{\mathrm{loc}}(\omega,\pi_1)$ лінійний простір усіх розподілів $u$
в області $\Omega$ таких, що $\chi u\in H^{s,s\gamma;\varphi}(\Omega)$ для кожної функції $\chi\in C^\infty (\overline\Omega)$ із $\mathrm{supp}\,\chi\subset\omega\cup\pi_1$.
Топологія у цьому просторі задається напівнормами
$$
u\mapsto\|\chi u\|_{H^{s,s\gamma;\varphi}(\Omega)},
$$
де $\chi$ -- довільна вище згадана функція.
Подібно до цього, позначимо через $H^{s,s\gamma;\varphi}_{\mathrm{loc}}(\pi_2)$ лінійний простір усіх розподілів $v$ на $S$ таких, що $\chi v\in H^{s,s\gamma;\varphi}(S)$ для кожної функції $\chi\in C^\infty (\overline S)$ із $\mathrm{supp}\,\chi\subset\pi_2$.
Топологія у цьому просторі задається напівнормами
$$
v\mapsto\|\chi v\|_{H^{s,s\gamma;\varphi}(S)},
$$
де $\chi$ -- довільна тільки згадана функція.
Нарешті, позначимо через $H^{s;\varphi}_{\mathrm{loc}}(\pi_3)$ лінійний простір усіх розподілів $w$ на $G$ таких, що $\chi w\in H^{s;\varphi}(G)$ для кожної функції $\chi\in C^\infty (\overline G)$ із $\mathrm{supp}\,\chi\subset\pi_3$.
Топологія у цьому просторі задається напівнормами
$$
w\mapsto\|\chi w\|_{H^{s;\varphi}(G)},
$$
де $\chi$ -- довільна тільки згадана функція.

\textbf{Теорема 2.} \it
Нехай
$u\in H^{\sigma_0,\sigma_0/(2b)}(\Omega)$ є узагальненим розв'язком параболічної
задачі \eqref{16f1}--\eqref{16f3} з правими частинами \eqref{16f12}.
Припустимо, що
\begin{align}\label{16f13}
f&\in H^{\sigma-2m,(\sigma-2m)/(2b);\varphi}_{\mathrm{loc}}
(\omega,\pi_1),
\\ \label{16f14}
g_{j}&\in H^{\sigma-m_j-1/2,(\sigma-m_j-1/2)/(2b);\varphi}_{\mathrm{loc}}
(\pi_2),\quad\mbox{для всіх}\quad j\in\{1,\dots,m\},
\\ \label{16f15}
h_{k}&\in H^{\sigma-2bk-b;\varphi}_{\mathrm{loc}}
(\pi_3),\quad\mbox{для всіх}\quad k\in\{0,\dots,\varkappa-1\},
\end{align}
для деяких $\sigma>\sigma_0$ і $\varphi\in\mathcal{M}$.
Тоді $u\in H^{\sigma,\sigma/(2b);\varphi}_{\mathrm{loc}}(\omega,\pi_1)$.
\rm

У випадку, коли $\pi_1=\varnothing$,  теорема~2
стверджує, що гладкість розв'язку підвищується в околах внутрішніх точок
замкненої області $\overline{\Omega}$.
Якщо $\pi_3=\varnothing$ то ця теорема стверджує про підвищення гладкості
розв'язку при $t>0$ і є наслідком теореми 4.3 з \cite{LosMikhMurach15Math}.
Якщо $\pi_1=\partial\Omega\setminus\Gamma$, $\pi_2=S$, $\pi_3=G$, то
підвищення гладкості розв'язку відбувається на множині $\overline\Omega\setminus\Gamma$.

У п.6 виведемо теорему~2 з теореми~1.

\textbf{4. Інтерполяція з функціональним параметром.}
Нагадаємо означення інтерпо\-ля\-ції з функціональним параметром у
випадку загальних гільбертових просторів та обговоримо інтерполяційні властивості,
які будуть використані у наступному розділі.
Ми слідуємо монографії
\cite[п.
1.1]{MikhailetsMurach14} (див. також \cite[п.~2]{08MFAT1}).
Обмежимось розглядом випадку сепарабельних комплексних гільбертових просторів.

Нехай $X:=[X_{0},X_{1}]$ є впорядкованою парою
сепарабельних комплексних гільбертових просторів для яких
має місце неперервне і щільне вкладання $X_{1}\hookrightarrow X_{0}$.
Таку пару називають допустимою.
Для неї існує оператор $J$, такий, що він є
самоспряженим додатньо визначеним оператором в
$X_{0}$ з областю визначення $X_{1}$, причому $\|Ju\|_{X_0}=\|u\|_{X_1}$.
Оператор $J$  визначається
парою $X$ однозначно і називається породжуючим оператором для~$X$.

Нехай $\psi\in\mathcal{B}$. Тут через $\mathcal{B}$ позначено множину всіх вимірних
за Борелем функцій $\psi:(0,\infty)\rightarrow(0,\infty)$, для яких $\psi$ є
обмеженою на кожному відрізку $[a,b]$, де $0<a<b<\infty$, і  $1/\psi$ є
обмеженою на кожному промені $[a,\infty)$, де $a>0$.

Розглянемо оператор $\psi(J)$, він є додатньо визначеним оператором в $X_{0}$ як
борелевська функція $\psi$ від $J$. Позначимо через $[X_{0},X_{1}]_{\psi}$, або скорочено
$X_{\psi}$, область візначення оператора $\psi(J)$, наділену скалярним добутком
$$
(u_{1},u_{2})_{X_{\psi}}:=(\psi(J)u_{1},\psi(J)u_{2})_{X_{0}}.\quad
$$
Він породжує норму $\|u\|_{X_{\psi}}:=\|\psi(J)u\|_{X_{0}}$. Простір $X_{\psi}$
є сепарабельним гільбертовим.

Функцію $\psi\in\mathcal{B}$ назвемо інтерполяційним параметром, якщо
для всіх допустимих пар $X=[X_{0},X_{1}]$ та $Y=[Y_{0},Y_{1}]$ гільбертових
просторів і для довільного лінійного відображення $T$, заданого на $X_{0}$
вірно таке: якщо звуження відображення $T$ на $X_{j}$ є обмеженим
оператором $T:X_{j}\rightarrow Y_{j}$ для кожного $j\in\{0,1\}$, тоді звуження
відображення $T$ на $X_{\psi}$ є також обмеженим оператором
$T:X_{\psi}\rightarrow Y_{\psi}$.

Якщо $\psi$ є інтерполяційним параметром, тоді будемо казати, що гільбертів простір
$X_{\psi}$ отримано в результаті інтерполяції з функціональним параметром $\psi$
пари $X=[X_{0},X_{1}]$. В цьому випадку маємо щільні та неперервні вкладання
$X_{1}\hookrightarrow X_{\psi}\hookrightarrow X_{0}$.

Відомо, що функція $\psi\in\mathcal{B}$ є інтерполяційним параметором тоді і тільки тоді
коли $\psi$ є псевдоугнутою в околі $\infty$, тобто коли існує угнута додатна
функція $\psi_{1}(r)$ при $r\gg1$ така, що обидві функції
$\psi/\psi_{1}$ та $\psi_{1}/\psi$ є обмеженими в деякому околі $\infty$.
Цей критерій випливає з опису Ж.~Петре класу всіх інтерполяційних функцій
для вагових просторів типу $L_{p}(\mathbb{R}^{n})$ (див. \cite[Теорема
5.4.4]{BerghLefstrem76}). Доведення цього критерію наведено, наприклад, в  \cite[п. 1.1.9]{MikhailetsMurach14}.

Для нас важливим є такий наслідок з цього критерію \cite[Теорема
1.11]{MikhailetsMurach14}.

\textbf{Твердження~1.} \it
Припустимо, що функція $\psi\in\mathcal{B}$ є правильно змінною на нескінченно\-сті
функцією порядку $\theta$, де  $0<\theta<1$, тобто
$$
\lim_{r\rightarrow\infty}\;\frac{\psi(\lambda r)}{\psi(r)}=
\lambda^{\theta}\quad\mbox{для кожного}\quad\lambda>0.
$$
Тоді $\psi$ є інтерполяційним параметром.
\rm

У випадку степеневих функцій це твердження приводить до класичного результату
Ж.-Л.~Ліонса та С.~Г.~Крейна, який полягає в тому, що функція $\psi(r)\equiv
r^{\theta}$ є інтерполяційним параметром при $0<\theta<1$. Тут показник
$\theta$ розглядається як числовий параметр інтерполяції
(ця інтерполяція була застосована у формулі \eqref{16f10}).

У кінці цього пункту сформулюємо три властивості інтерполяції, які будуть використані
у подальших доведеннях (див. \cite[пп. 1.1.3, 1.1.5, 1.1.6]{MikhailetsMurach14}
та \cite[п. 1.17]{Triebel95}). Перша з них дозволяє звести інтерполяцію підпросторів або
фактор--просторів до інтерполяції вихідних просторів.
Зазначимо, що підпростори припускаються замкнені, і розглядаються, взагалі кажучи, не ортогональні проектори на підпростори.

\textbf{Твердження~2.} \it
Нехай $X=[X_{0},X_{1}]$ є допустимою парою гільбертових просторів, а $Y_{0}$ є
підпростором в $X_{0}$. Тоді $Y_{1}:=X_{1}\cap Y_{0}$ є підпростором в $X_{1}$.
Припустимо, що існує лінійне відображення $P:X_{0}\rightarrow X_{0}$, яке для
кожного $j\in\{0,\,1\}$ є проектором простору $X_{j}$ на його підпростір $Y_{j}$.
Тоді пари $[Y_{0},Y_{1}]$ та $[X_{0}/Y_{0},X_{1}/Y_{1}]$
є допустимими і для довільного інтерполяційного параметра $\psi\in\mathcal{B}$
правильні рівності
\begin{gather*}
[Y_{0},Y_{1}]_{\psi}=X_{\psi}\cap Y_{0},\\
[X_{0}/Y_{0},X_{1}/Y_{1}]_{\psi}=X_{\psi}/(X_{\psi}\cap Y_{0})
\end{gather*}
з еквівалентністю норм.
Тут $X_{\psi}\cap Y_{0}$ є підпростором у $X_{\psi}$.
\rm

Друга властивість дозволяє звести інтерполяцію прямих сум гільбертових просторів
до інтерполяції їх доданків.

\textbf{Твердження~3.} \it
Нехай $[X_{0}^{(j)},X_{1}^{(j)}]$, де $j=1,\ldots,p$, є скінченний набір допустимих пар гільбертових просторів. Тоді
$$
\biggl[\,\bigoplus_{j=1}^{p}X_{0}^{(j)},\,\bigoplus_{j=1}^{p}X_{1}^{(j)}\biggr]_{\psi}=\,
\bigoplus_{j=1}^{p}\bigl[X_{0}^{(j)},\,X_{1}^{(j)}\bigr]_{\psi}
$$
з рівністю норм. Тут функція $\psi\in\mathcal{B}$ є довільним інтерполяційним параметром.
\rm

Третя показує, що повторне застосування інтерполяції з функціональним параметром
дає знову інтерполяцію з деяким функціональним параметром.

\textbf{Твердження~4.} \it Нехай $\alpha,\beta,\psi\in\mathcal{B}$ і функція
$\alpha/\beta$ є обмеженою в околі нескінченності. Тоді для кожної допустимої пари $X$ гільбертових просторів правильна рівність
$[X_{\alpha},X_{\beta}]_{\psi}=X_{\omega}$ з рівністю норм.
Тут функція $\omega\in\mathcal{B}$ визначена за формулою
$\omega(t):=\alpha(t)\psi(\beta(t)/\alpha(t))$ при $t>0$.
Більше того, якщо $\alpha,\beta,\psi$ є інтерполяційними параметрами, то
$\omega$ також є інтерполяційним параметром.
\rm

\textbf{5. Інтерполяція пар гільбертових просторів Соболєва.}
У цьому пункті доведемо, що простори Хермандера, які з'являються в теоремі~1, можна
отримати шляхом інтерполяції з функціональним параметром їх соболєвських аналогів.

Нехай
\begin{equation}\label{16f16}
s,s_{0},s_{1},\gamma\in\mathbb{R},\quad s_{0}<s<s_{1},\quad\gamma>0,
\quad\mbox{і}\quad\varphi\in\mathcal{M}.
\end{equation}
Розглянемо функцію
\begin{equation}\label{16f17}
\psi(r):=
\begin{cases}
\;r^{(s-s_{0})/(s_{1}-s_{0})}\,\varphi(r^{1/(s_{1}-s_{0})})&\text{для}\quad r\geq1, \\
\;\varphi(1) & \text{для}\quad0<r<1.
\end{cases}
\end{equation}
За твердженням~1, ця функція  є інтерполяційним параметром, оскільки
вона є правильно змінною функцією на нескінченності порядку $\theta:=(s-s_{0})/(s_{1}-s_{0})$
з $0<\theta<1$. Надалі будемо інтерполювати пари соболєвських просторів з
функціональним параметром $\psi$.

Попередньо запишемо необхідні інтерполяційні формули для просторів Хермандера,
заданих на $\mathbb{R}^{k}$, де ціле $k\geq1$, $G$, $\Gamma$ та $S$. Ці формули
було доведено в інших роботах.

Інтерполяція ізотропних соболєвських просторів досліджена в роботах
В.~А.~Ми\-хайлеця  та О.О. Мурача \cite{12BJMA2,MikhailetsMurach14}.
У припущенні \eqref{16f16} правильні такі рівності \cite[теореми 1.14, 2.2, 3.2]{MikhailetsMurach14}:
\begin{align}\label{16f18}
&H^{s;\varphi}(\mathbb{R}^{k})=\bigl[H^{s_{0}}(\mathbb{R}^{k}),
H^{s_{1}}(\mathbb{R}^{k})\bigr]_{\psi},\\
\label{16f19}
&H^{s;\varphi}(G)=\bigl[H^{s_{0}}(G),
H^{s_{1}}(G)\bigr]_{\psi},\\
\label{16f20}
&H^{s;\varphi}(\Gamma)=\bigl[H^{s_{0}}(\Gamma),
H^{s_{1}}(\Gamma)\bigr]_{\psi}
\end{align}
з точністю до еквівалентності норм.

Подібна рівність виконується  для анізотропних просторів $H^{s,s\gamma;\varphi}(S)$.
Додатково до \eqref{16f16} припустимо, що $s_0\geq0$. Тоді правильна така рівність
\cite[теорема~2]{Los15NK2}:
\begin{equation}\label{16f21}
H^{s,s\gamma;\varphi}(S)=\bigl[H^{s_{0},s_{0}\gamma}(S),
H^{s_{1},s_{1}\gamma}(S)\bigr]_{\psi}
\end{equation}
з точністю до еквівалентності норм.

Встановимо версію цих формул для простору $H^{s,s\gamma;\varphi}(\Omega)$.

\textbf{Лема 1.} \it Припустимо додатково до \eqref{16f16}, що $s_0\geq0$.
Тоді правильна така інтерполяційна формула
\begin{equation}\label{16f22}
H^{s,s\gamma;\varphi}(\Omega)=\bigl[H^{s_{0},s_{0}\gamma}(\Omega),
H^{s_{1},s_{1}\gamma}(\Omega)\bigr]_{\psi}
\end{equation}
з точністю до еквівалентності норм.
\rm

\textbf{Доведення.}
Доведення формули \eqref{16f22} спирається на твердження~2 у частині інтерполяції
фактор--просторів і на інтерполяційну формулу
(див. \cite[Лема 1]{Los15NK2} або \cite[Лема 7.1]{LosMikhMurach15Math})
\begin{equation}\label{16f23}
H^{s,s\gamma;\varphi}(\mathbb{R}^{k+1})=\bigl[H^{s_{0},s_{0}\gamma}(\mathbb{R}^{k+1}),
H^{s_{1},s_{1}\gamma}(\mathbb{R}^{k+1})\bigr]_{\psi},
\end{equation}
яка виконується з рівністю норм.
Тут $H^{s,s\gamma;\varphi}(\mathbb{R}^{k+1}):=H^{\mu}(\mathbb{R}^{k+1})$,
де показник $\mu$ визначений формулою \eqref{16f70}, у якій $k:=k+1$.

Як зазначалося у п. 2, простір
$H^{s,s\gamma;\varphi}(\Omega)$
є факторпростором простору
$H^{s,s\gamma;\varphi}(\mathbb{R}^{n+1})$ за його підпростором
$H^{s,s\gamma;\varphi}(\mathbb{R}^{n+1},\Omega)$ \eqref{16f69}.
Для виведення \eqref{16f22} з \eqref{16f23}
побудуємо проектор $P$ кожного простору
$H^{s_{j},s_{j}\gamma}(\mathbb{R}^{n+1})$, із $j\in\{0,1\}$,
на його підпростір $H^{s_{j},s_{j}\gamma}_{Q}(\mathbb{R}^{n+1})$.

Як відомо \cite[п.~9.6]{BesovIlinNikolskii75}, існує лінійний обмежений оператор продовження $T_{\Omega}\,:\,L_2(\Omega)\rightarrow  L_2(\mathbb{R}^{n+1})$ такий, що звуження $T$
на кожний простір $H^{\sigma,\sigma\gamma}(\Omega)$,
де $\sigma$ і $\sigma\gamma$ натуральні числа, є обмеженим оператором
$T_{\Omega}\,:\,H^{\sigma,\sigma\gamma}(\Omega)\rightarrow
H^{\sigma,\sigma\gamma}(\mathbb{R}^{n+1}).$
Виходячи з цього припустимо додатково, що всі $s_{0}$, $s_{1}$,
$s_{0}\gamma$, та $s_{1}\gamma$ є цілими невід'ємними числами.
Розглянемо відображення $P :\,w\mapsto w-T_\Omega(w\!\!\upharpoonright\!\Omega)$, де $w\in H^{s_{0},s_{0}\gamma}(\mathbb{R}^{n+1})$.
Оператор $P$ є проектором.
Справді, $P$ є лінійним обмеженим оператором на
$H^{s_{j},s_{j}\gamma}(\mathbb{R}^{n+1})$ для кожного $j\in\{0,1\}$.
Більш того, якщо $w=0$ в $\Omega$, то
$T_\Omega(w\!\!\upharpoonright\!\Omega)=0$ в $\mathbb{R}^{n+1}$;
тому $Pw=w$ для кожної
$w\in H^{s_{j},s_{j}\gamma}_{Q}(\mathbb{R}^{n+1})$.

Оскільки проектор $P$ заданий, можемо використати твердження~2
та формулу \eqref{16f23} з $k:=n$ і записати
\begin{align*}
\bigl[&H^{s_{0},s_{0}\gamma}(\Omega),
H^{s_{1},s_{1}\gamma}(\Omega)\bigr]_{\psi}\\
&=\bigl[H^{s_{0},s_{0}\gamma}(\mathbb{R}^{n+1})/
H^{s_{0},s_{0}\gamma}_{Q}(\mathbb{R}^{n+1}),
H^{s_{1},s_{1}\gamma}(\mathbb{R}^{n+1})/
H^{s_{1},s_{1}\gamma}_{Q}(\mathbb{R}^{n+1})\bigr]_{\psi}\\
&=\bigl[H^{s_{0},s_{0}\gamma}(\mathbb{R}^{n+1}),
H^{s_{1},s_{1}\gamma}(\mathbb{R}^{n+1})\bigr]_{\psi}\big/
\bigl(\bigl[H^{s_{0},s_{0}\gamma}(\mathbb{R}^{n+1}),
H^{s_{1},s_{1}\gamma}(\mathbb{R}^{n+1})\bigr]_{\psi}\cap
H^{s_{0},s_{0}\gamma}_{Q}(\mathbb{R}^{n+1})\bigr)\\
&=H^{s,s\gamma;\varphi}(\mathbb{R}^{n+1})/
\bigl(H^{s,s\gamma;\varphi}(\mathbb{R}^{n+1})\cap
H^{s_{0},s_{0}\gamma}_{Q}(\mathbb{R}^{n+1})\bigr)=
H^{s,s\gamma;\varphi}(\mathbb{R}^{n+1})/H^{s,s\gamma;\varphi}_{Q}(\mathbb{R}^{n+1})\\
&=H^{s,s\gamma;\varphi}(\Omega)
\end{align*}
з точністю до еквівалентності норм.

Формулу \eqref{16f22} встановлено у додатковому припущенні, що
всі $s_{0}$, $s_{1}$, $s_{0}\gamma$, та $s_{1}\gamma$ є цілими
невід'ємними числами. Це припущення можна прибрати, оскільки звуження
оператора $T_{\Omega}$ на кожний простір $H^{\sigma,\sigma\gamma}(\Omega)$,
де дійсне $\sigma\geq 0$, є обмеженим оператором
\begin{equation}\label{16f24}
T_{\Omega}\,:\,H^{\sigma,\sigma\gamma}(\Omega)\rightarrow
H^{\sigma,\sigma\gamma}(\mathbb{R}^{n+1}).
\end{equation}
Доведемо обмеженість оператора \eqref{16f24} за допомогою інтерполяції.
Нехай число $k_1\in\mathbb{N}$ таке, що $k_1>\sigma$ і $k_1\gamma\in\mathbb{N}$.
Розглянемо обмежені оператори
$$
T_{\Omega}\,:\,L_2(\Omega)\rightarrow  L_2(\mathbb{R}^{n+1})
\quad\mbox{і}\quad
T_{\Omega}\,:\,H^{k_1,k_1\gamma}(\Omega)\rightarrow  H^{k_1,k_1\gamma}(\mathbb{R}^{n+1}).
$$
Оскільки $\psi$ є інтерполяційним параметром, то з обмеженості останніх двох операторів випливає обмеженість
оператора
\begin{equation}\label{16f25}
T_{\Omega} :\ \left[L_2(\Omega),H^{k_1,k_1\gamma}(\Omega)\right]_{\psi}
\rightarrow\left[L_2(\mathbb{R}^{n+1}),H^{k_1,k_1\gamma}(\mathbb{R}^{n+1})\right]_{\psi}.
\end{equation}
Тепер обмеженість оператора \eqref{16f24} є прямим наслідком обмеженості
оператора \eqref{16f25} та інтерполяційних рівностей  \eqref{16f22} і \eqref{16f23}
при $\varphi=1$, $s_0=0$, $s=\sigma$, $s_1=k_1$, $k=n$.

Для завершення доведення формули \eqref{16f22} достатньо повторити частину її
доведення починаючи з побудови проектора $P$, скориставшись
при цьому оператором продов\-ження~\eqref{16f24}.

Лема 1 доведена.

Для доведення інтерполяційної формули для просторів
$\mathcal{Q}^{s-2m,(s-2m)/(2b);\varphi}$,
до яких належать праві частини розглядуваної задачі,
нам будуть потрібні такі дві леми про властивості
оператора даних Коші в анізотропних просторах Хермандера.

\textbf{Лема 2.} \it Нехай задане довільне $r\in\mathbb{N}$.
Для довільних  дійсного числа $s>2br-b$ і функції $\varphi\in\mathcal{M}$
лінійні відображення
\begin{equation}\label{16f26}
w\,\mapsto \, (w\!\upharpoonright_{t=0},\,
(\partial_t w)\!\upharpoonright_{t=0},\dots,
(\partial^{\,r-1}_t w)\!\upharpoonright_{t=0}),\quad\mbox{де}
\quad w\in C_0^{\infty}(\mathbb{R}^{n}),
\end{equation}
та
\begin{equation}\label{16f27}
g\,\mapsto \, (g\upharpoonright\!\Gamma,
(\partial_t g)\upharpoonright\!\Gamma,\dots,
(\partial^{\,r-1}_t g)\upharpoonright\!\Gamma),
\quad\mbox{де}\quad g\in C^{\infty}(\overline{S})
\end{equation}
продовжуються по неперервності до обмежених операторів
\begin{equation}\label{16f28}
R_{\mathbb{R}^{n-1}}\,:\, H^{s,s/(2b);\varphi}(\mathbb{R}^{n})
\rightarrow \bigoplus_{k=0}^{r-1}H^{s-2bk-b;\varphi}(\mathbb{R}^{n-1})
\end{equation}
та
\begin{equation}\label{16f29}
R_{\Gamma}\,:\, H^{s,s/(2b);\varphi}(S)\rightarrow \bigoplus_{k=0}^{r-1}
H^{s-2bk-b;\varphi}(\Gamma)
\end{equation}
відповідно.
\rm

\textbf{Доведення.}
У випадку $\varphi\equiv 1$ (простори Соболєва) тверждення леми є відомим
\cite[п.II, теореми 6, 7]{Slobodetskii58}
(див. також \cite[теорема 2.4]{ZhitarashuEidelman98}).
Звідси загальний випадок $\varphi\in\mathcal{M}$
виведемо за допомогою інтерполяції.
Зробимо це для відображення \eqref{16f26}.
Виберемо $s_0$ і $s_1$ так, щоб $2br-b<s_0<s<s_1$.
Маємо лінійні обмежені оператори
$$
R_{\mathbb{R}^{n-1}}\,:\, H^{s_j,s_j/(2b)}(\mathbb{R}^{n})\rightarrow
\bigoplus_{k=0}^{r-1}H^{s_j-2bk-b}(\mathbb{R}^{n-1})
\quad\mbox{для всіх}\quad j\in\{0,1\}.
$$
Застосувавши тут інтерполяцію з функціональним параметром
$\psi$, ми отримаємо обмежений оператор
\begin{equation}\label{16f30}
\begin{split}
R_{\mathbb{R}^{n-1}}\,:\, \bigl[H^{s_0,s_0/(2b)}(\mathbb{R}^{n}),&
H^{s_1,s_1/(2b)}(\mathbb{R}^{n})\bigr]_{\psi}\rightarrow \\ &
\biggl[\,\bigoplus_{k=0}^{r-1}H^{s_0-2bk-b}(\mathbb{R}^{n-1}),
\bigoplus_{k=0}^{r-1}H^{s_1-2bk-b}(\mathbb{R}^{n-1})\,\biggr]_{\psi}.
\end{split}
\end{equation}
Оператор \eqref{16f30} є розширенням за неперервністю відображення \eqref{16f26}
оскільки множина $C_0^{\infty}(\mathbb{R}^{n})$ щільна в області визначення оператора \eqref{16f30}.
З обмеженості оператора \eqref{16f30}, твердження 3 та рівностей  \eqref{16f18} і \eqref{16f23}
при $k=n-1$ випливає обмеженість оператора \eqref{16f28}.

Для відображення \eqref{16f27} доведення
проводиться за тією ж схемою з використанням інтерполяційних
формул \eqref{16f20} і \eqref{16f21}.

Лема~2 доведена.

\textbf{Лема 3.} \it Нехай задано довільне $r\in\mathbb{N}$.
Тоді існує лінійне відображення $T:\,(L_2(\Gamma))^r\to
L_2(S)$ таке, що для довільних дійсного числа $s>2br-b$ і функціонального параметра
$\varphi\in\mathcal{M}$ звуження відображення $T$ на простір
$$
\bigoplus_{k=0}^{r-1}H^{s-2bk-b;\varphi}(\Gamma)
$$
є лінійним обмеженим оператором
\begin{equation}\label{16f31}
T:\bigoplus_{k=0}^{r-1}H^{s-2bk-b;\varphi}(\Gamma)
\rightarrow H^{s,s/(2b);\varphi}(S),
\end{equation}
правим оберненим до оператора \eqref{16f29}.
\rm

\textbf{Доведення.} Спочатку доведемо цю лему
для відповідних гільбертових просторів, заданих на
$\mathbb{R}^{n-1}$ та $\mathbb{R}^{n}$ замість $\Gamma$ та $S$.
За аналогією з \cite{Hermander63} (доведення теореми 2.5.7)
розглянемо лінійне відображення
\begin{equation*}
T_1\,:\,v\mapsto F_{\xi}^{-1}\bigl(
\beta(\langle\xi\rangle^{2b}t)\,\sum_{k=0}^{r-1}
\frac{1}{k!}\,\widehat{v_k}(\xi)t^k\bigr),
\end{equation*}
задане на вектор-функціях
$$
v:=(v_0,\dots,v_{r-1})\in\bigl(S'(\mathbb{R}^{n-1})\bigr)^r.
$$
Тут $\beta\in C^{\infty}_{0}(\mathbb{R})$ є деяка функція, що дорівнює $1$ в околі нуля,
$\langle\xi\rangle:=(1+|\xi|^2)^{1/2}$ -- згладжений модуль вектора $\xi\in\mathbb{R}^{n-1}$,
$\widehat{w}(\xi):=(Fw)(\xi)$, де $\xi\in\mathbb{R}^{n-1}$, є (повне) перетворення
Фур'є розподілу $w\in S'(\mathbb{R}^{n-1})$, а $F^{-1}_\xi$~--~оператор оберненого
перетворення Фур'є по частотній змінній $\xi\in\mathbb{R}^{n-1}$.
Звісно, відображення $T_1$
є лінійним неперервним оператором
\begin{equation}\label{16f32}
 T_1\,:\,\bigl(S'(\mathbb{R}^{n-1})\bigr)^r\rightarrow
 S'(\mathbb{R}^{n}).
\end{equation}
Покажемо, що звуження оператора \eqref{16f32} на простір
$$
\bigoplus_{k=0}^{r-1}H^{s-2bk-b;\varphi}(\mathbb{R}^{n-1})
$$
є обмеженим оператором
\begin{equation}\label{16f33}
T_1\,:\, \bigoplus_{k=0}^{r-1}H^{s-2bk-b;\varphi}(\mathbb{R}^{n-1})
\rightarrow H^{s,s/(2b);\varphi}(\mathbb{R}^{n}),
\end{equation}
правим оберненим до оператора \eqref{16f28}.

Попередньо відмітимо, що
\begin{equation}\label{16f34}
R_{\mathbb{R}^{n-1}}T_1v=v \quad\mbox{для довільних}\quad
v\in\bigl(S(\mathbb{R}^{n-1})\bigr)^r.
\end{equation}
Справді, для кожного індекса $j\in\{0,\dots,r-1\}$ і
довільної вектор-функції
$$
v=(v_0,\dots,v_{r-1})\in\bigl(S(\mathbb{R}^{n-1})\bigr)^r
$$
виконуються рівності
\begin{align*}
\widehat{(\partial^j_tT_1v)\!\upharpoonright_{t=0}}&=
\partial^j_t(\widehat{T_1v})\!\upharpoonright_{t=0}=
\partial^j_t\biggl(
\beta(\langle\xi\rangle^{2b}t)\,\sum_{k=0}^{r-1}
\frac{1}{k!}\,\widehat{v_k}(\xi)t^k\biggr)\!\upharpoonright_{t=0}\\
&=\beta(0)\biggl(\partial^j_t\bigl(\sum_{k=0}^{r-1}
\frac{1}{k!}\,\widehat{v_k}(\xi)t^k\bigr)
\biggr)\!\upharpoonright_{t=0}=\beta(0)\,j!\,\frac{1}{j!}\,\widehat{v_j}(\xi)=
\widehat{v_j}(\xi).
\end{align*}
Тут враховано, що $\beta=1$ в околі нуля.
Отже, $\widehat{R_{\mathbb{R}^{n-1}}T_1v}=\widehat{v}$, що
рівносильно \eqref{16f34}.

Доведемо спочатку обмеженість оператора \eqref{16f33} у соболєвському випадку,
коли $s=2bm$, де $m\in\mathbb{N}\cup\{0\}$, і $\varphi\equiv1$.
Відповідну оцінку встановимо використавши у просторі $H^{2bm,m}(\mathbb{R}^{n})$
таку еквівалентну норму \cite[п.9.1]{BesovIlinNikolskii75}:
\begin{equation*}
||u||^2_{2bm,m}:=||u||^2_{L_2(\mathbb{R}^{n})}+
\sum_{j=1}^{n-1}||\partial_{x_j}^{2bm}u||^2_{L_2(\mathbb{R}^{n})}+
||\partial_{t}^{m}u||^2_{L_2(\mathbb{R}^{n})}.
\end{equation*}
Для довільного $v\in\bigl(S(\mathbb{R}^{n-1})\bigr)^r$ запишемо
на підставі рівності Парсеваля
\begin{align*}
||T_1v||^2_{2bm,m}&=
||T_1v||^2_{L_2(\mathbb{R}^{n})}+
\sum_{j=1}^{n-1}||\partial_{x_j}^{2bm}(T_1v)||^2_{L_2(\mathbb{R}^{n})}+
||\partial_{t}^{m}(T_1v)||^2_{L_2(\mathbb{R}^{n})}\\
&=||\widehat{T_1v}||^2_{L_2(\mathbb{R}^{n})}+
\sum_{j=1}^{n-1}||\xi_j^{2bm}\widehat{T_1v}||^2_{L_2(\mathbb{R}^{n})}+
||\partial_{t}^{m}(\widehat{T_1v})||^2_{L_2(\mathbb{R}^{n})}\\
&\leq\sum_{k=0}^{r-1}\frac{1}{k!}\,
||\beta(\langle\xi\rangle^{2b}t)\,\widehat{v_k}(\xi)t^k||^2_{L_2(\mathbb{R}^{n})}
+\sum_{j=1}^{n-1}
\sum_{k=0}^{r-1}\frac{1}{k!}\,
||\xi_j^{2bm}\beta(\langle\xi\rangle^{2b}t)\,\widehat{v_k}(\xi)t^k||^2_{L_2(\mathbb{R}^{n})}\\
&\quad+\sum_{k=0}^{r-1}\frac{1}{k!}\,||\partial_{t}^{m}(
\beta(\langle\xi\rangle^{2b}t)\,\widehat{v_k}(\xi)t^k)||^2_{L_2(\mathbb{R}^{n})}.
\end{align*}
Оцінимо окремо кожну з останніх трьох норм. Використавши заміну змінних
$\tau=\langle\xi\rangle^{2b}t$, де $\xi$ фіксовано, отримуємо рівності
\begin{align*}
||\partial_{t}^{m}(\beta(\langle\xi\rangle^{2b}t)\widehat{v_k}(\xi)t^k)||^2_{L_2(\mathbb{R}^{n})}&=
\int\limits_{\mathbb{R}^{n}}
|\widehat{v_k}(\xi)|^2|\partial_{t}^{m}(\beta(\langle\xi\rangle^{2b}t)t^k)|^2d\xi dt\\
&=\int\limits_{\mathbb{R}^{n-1}}
|\widehat{v_k}(\xi)|^2d\xi \int\limits_{\mathbb{R}}
|\partial_{t}^{m}(\beta(\langle\xi\rangle^{2b}t)t^k)|^2 dt\\
&=\int\limits_{\mathbb{R}^{n-1}}
|\widehat{v_k}(\xi)|^2\langle\xi\rangle^{4bm-4bk-2b} d\xi \int\limits_{\mathbb{R}}
|\partial_{\tau}^{m}(\beta(\tau)\tau^{k})|^2 d\tau.
\end{align*}
Тобто
$$
||\partial_{t}^{m}(\beta(\langle\xi\rangle^{2b}t)\widehat{v_k}(\xi)t^k)||^2_{L_2(\mathbb{R}^{n})}=
c_1\,||v_k||^2_{H^{2bm-2bk-b}(\mathbb{R}^{n-1})},
$$
де
$$
c_1:=\int\limits_{\mathbb{R}}
|\partial_{\tau}^{m}(\beta(\tau)\tau^{k})|^2 d\tau<\infty.
$$
Далі
\begin{align*}
||\xi_j^{2bm}\beta(\langle\xi\rangle^{2b}t)\widehat{v_k}(\xi)t^k||^2_{L_2(\mathbb{R}^{n})}&=
\int\limits_{\mathbb{R}^{n}}
|\xi_j|^{4bm}|\beta(\langle\xi\rangle^{2b}t)|^2|\widehat{v_k}(\xi)|^2|t|^{2k}d\xi dt\\
&=\int\limits_{\mathbb{R}^{n-1}}
|\xi_j|^{4bm}|\widehat{v_k}(\xi)|^2d\xi \int_{\mathbb{R}}
|\beta(\langle\xi\rangle^{2b}t)|^2|t|^{2k} dt\\
&=\int\limits_{\mathbb{R}^{n-1}}
|\xi_j|^{4bm}|\widehat{v_k}(\xi)|^2\langle\xi\rangle^{-4bk-2b} d\xi \int_{\mathbb{R}}
|\beta(\tau)|^2|\tau|^{2k} d\tau\\
&=c_2\,\int\limits_{\mathbb{R}^{n-1}}
|\widehat{v_k}(\xi)|^2\langle\xi\rangle^{4bm-4bk-2b}
|\xi_j|^{4bm}\langle\xi\rangle^{-4bm}d\xi\\
&\leq c_2\,\int\limits_{\mathbb{R}^{n-1}}
|\widehat{v_k}(\xi)|^2\langle\xi\rangle^{4bm-4bk-2b}d\xi.
\end{align*}
Отже,
$$
||\xi_j^{2bm}\beta(\langle\xi\rangle^{2b}t)\widehat{v_k}(\xi)t^k||^2_{L_2(\mathbb{R}^{n})}
\leq c_2\,||v_k||^2_{H^{2bm-2bk-b}(\mathbb{R}^{n-1})},
$$
де
$$
c_2:=\int_{\mathbb{R}}
|\beta(\tau)|^2|\tau|^{2k} d\tau<\infty.
$$
Нарешті, замінивши у попередній оцінці $\xi_j^{2bm}$ на число $1$,
дістаємо оцінку для першої норми
$$
||\beta(\langle\xi\rangle^{2b}t)\widehat{v_k}(\xi)t^k||^2_{L_2(\mathbb{R}^{n})}
\leq c_2\,||v_k||^2_{H^{2bm-2bk-b}(\mathbb{R}^{n-1})}.
$$

Таким чином,
\begin{equation*}
||T_1v||_{H^{2bm,m}(\mathbb{R}^{n})}\leq c\,\sum_{k=0}^{r-1}
||v_k||^2_{H^{2bm-2bk-b}(\mathbb{R}^{n-1})}
\end{equation*}
для довільних $v\in\bigl(S(\mathbb{R}^{n-1})\bigr)^r$
з деякою сталою $c>0$, що не залежить від $v$.
Оскільки множина $\bigl(S(\mathbb{R}^{n-1})\bigr)^r$ є щільною
в просторі
$$
\bigoplus_{k=0}^{r-1}H^{2bm-2bk-b}(\mathbb{R}^{n-1}),
$$
то цим і доведено обмеженість оператора
\begin{equation}\label{16f35}
T_1\,:\, \bigoplus_{k=0}^{r-1}H^{2bm-2bk-b}(\mathbb{R}^{n-1})
\rightarrow H^{2bm,m}(\mathbb{R}^{n}).
\end{equation}

Звідси виведемо обмеженість оператора \eqref{16f33} за допомогою інтерполяції
з функ\-ціо\-нальним параметром відповідних пар просторів Соболєва.

Нехай $s_0=0$, а кратне $2b$ число $s_1$ таке, що $s_1>s$.
На підставі \eqref{16f35} маємо лінійні обмежені оператори
\begin{equation}\label{16f36}
T_1\,:\,\bigoplus_{k=0}^{r-1}H^{s_j-2bk-b}(\mathbb{R}^{n-1})
\,\rightarrow\,  H^{s_j,s_j/(2b)}(\mathbb{R}^{n})
\quad\mbox{для кожного}\quad j\in\{0,1\}.
\end{equation}
Застосувавши тут інтерполяцію з функціональним параметром
$\psi$ і скориставшись твердженням 3, отримаємо обмежений оператор
\begin{equation}\label{16f37}
T_1\,:\,\bigoplus_{k=0}^{r-1}
\bigl[H^{s_0-2bk-b}(\mathbb{R}^{n-1}),H^{s_1-2bk-b}(\mathbb{R}^{n-1})\bigr]_{\psi}
\,\rightarrow\,\bigl[H^{s_0,s_0/(2b)}(\mathbb{R}^{n}),H^{s_1,s_1/(2b)}(\mathbb{R}^{n})\bigr]_{\psi}.
\end{equation}
Він є звуженням оператора \eqref{16f36} при $j=0$.

В результаті обмеженість оператора \eqref{16f33} є прямим наслідком обмеженості
оператора \eqref{16f37} і інтерполяційних рівностей  \eqref{16f18} та \eqref{16f23}.

Тепер рівність \eqref{16f34} продовжується за неперервністю на всі вектори
$$
v\in\bigoplus_{k=0}^{r-1}H^{s-2bk-b;\varphi}(\mathbb{R}^{n-1}),
$$
тобто оператор \eqref{16f33} є правим оберненим  до
оператора \eqref{16f28}.

Для подальших міркувань нам будуть потрібні такі аналоги операторів
\eqref{16f28} і \eqref{16f33}.

Для кожної функції $w\in H^{s,s/(2b);\varphi}(\Pi)$ покладемо
$R_{\Pi}w:=R_{\mathbb{R}^{n-1}}u$, де функція
$u\in H^{s,s/(2b);\varphi}(\mathbb{R}^{n})$ така, що
$u\upharpoonright_{\Pi}=w$. Це означення коректне, тобто не
залежить від зробленого вибору функції $u$. Лінійне відображення
$w\,\mapsto\,R_{\Pi}w$ задає обмежений оператор
\begin{equation}\label{16f38}
R_{\Pi}\,:\, H^{s,s/(2b);\varphi}(\Pi)
\rightarrow \bigoplus_{k=0}^{r-1}H^{s-2bk-b;\varphi}(\mathbb{R}^{n-1}).
\end{equation}
Це є прямим наслідком обмеженості оператора \eqref{16f28}
і означення норми у просторі $H^{s,s/(2b);\varphi}(\Pi)$.

Для кожного вектора
$$
v:=(v_0,v_1,\dots,v_{r-1})\in
\bigoplus_{k=0}^{r-1}H^{s-2bk-b;\varphi}(\mathbb{R}^{n-1})
$$
покладемо $T_{\Pi}v:=\bigl(T_1v\bigr)\upharpoonright_{\Pi}$.
Звісно, лінійне відображення $v\mapsto T_{\Pi}v$ є обмеженим оператором
\begin{equation}\label{16f39}
T_\Pi\,:\,
\bigoplus_{k=0}^{r-1}H^{s-2bk-b;\varphi}(\mathbb{R}^{n-1})
\,\rightarrow\, H^{s,s/(2b);\varphi}(\Pi).
\end{equation}
Відмітимо, що
$$
R_\Pi T_\Pi v=R_\Pi
\bigl(T_1v\bigr)\upharpoonright_{\Pi}=
R_{\mathbb{R}^{n-1}} T_1v=v
$$
для довільного
$$
v\in\bigoplus_{k=0}^{r-1}H^{s-2bk-b;\varphi}(\mathbb{R}^{n-1}).
$$
Отож оператор \eqref{16f39} є правим оберненим до
оператора \eqref{16f38}.

Тепер потрібний нам оператор \eqref{16f31} побудуємо
за допомогою оператора \eqref{16f39} та локального "розпрямлення"\
і "склеювання"\ многовидів $\Gamma$ та $S$.

Довільно виберемо скінченний атлас із
$C^{\infty}$-структури на $\Gamma$, утворений локальними
картами $\nobreak{\theta_{j}:\mathbb{R}^{n-1}\leftrightarrow \Gamma_{j}}$, де
$j=1,\ldots,\lambda$. Тут відкрити множини $\Gamma_{1},\ldots,\Gamma_{\lambda}$
складають покриття многовиду $\Gamma$. Окрім цього, довільно виберемо функції
$\chi_{j}\in C^{\infty}(\Gamma)$, $j=1,\ldots,\lambda$, такі, що
$\mathrm{supp}\,\chi_{j}\subset\Gamma_{j}$ і $\sum_{j=1}^{\lambda}\chi_{j}\equiv1$
на $\Gamma$.

За означенням простору $H^{\sigma;\varphi}(\Gamma)$,
лінійне відображення "розпрямлення"
\begin{equation*}
L:\,v\mapsto(\,(\chi_{1}v)\circ\theta_{1},\ldots,(\chi_{\lambda}v)\circ\theta_{\lambda}\,),
\end{equation*}
де $v\in L_2(\Gamma)$, задає ізометричний оператор
\begin{equation}\label{16f40}
L:\,H^{\sigma;\varphi}(\Gamma)\rightarrow
\bigl(H^{\sigma;\varphi}(\mathbb{R}^{n-1})\bigr)^{\lambda}
\end{equation}
для кожного $\sigma>0$. Тут і далі, як звичайно, символ
"$\circ$"\ позначає композицію відображень.

Розглянемо лінійне відображення "склеювання"
$$
K:\,(h_{1},\ldots,h_{\lambda})\mapsto\sum_{j=1}^{\lambda}\,
\Theta_{j}\left((\eta_{j}h_{j})\circ\theta_{j}^{-1}\right),
$$
де $h_{1},\ldots,h_{\lambda}\in L_2(\mathbb{R}^{n-1})$.
Тут для кожного номера $j=1,\ldots,\lambda$ функція $\eta_{j}\in
C_{0}^{\infty}(\mathbb{R}^{n-1})$ така, що $\eta_{j}=1$ на
множині $\theta^{-1}_{j}(\mathrm{supp}\,\chi_{j})$.
Тут також $\Theta_{j}$
оператор продовження нулем з $\Gamma_j$ на многовид $\Gamma$
(оператор $\Theta_{j}$ коректно означений на функціях,
носій яких лежить в $\Gamma_j$).
Це відображення задає обмежений оператор
\begin{equation*}
K :\,
\bigl(H^{\sigma;\varphi}(\mathbb{R}^{n-1})\bigr)^{\lambda}
\rightarrow H^{\sigma;\varphi}(\Gamma),
\end{equation*}
для кожного $\sigma>0$,
який є лівим оберненим до оператора \eqref{16f40}
(див. \cite{MikhailetsMurach14}, доведення теореми~2.2).

Для бічної поверхні $S$ аналогом $K$ служить таке лінійне
відображення:
\begin{equation*}
K_S\,:\,(g_1(\cdot,t),\dots,g_\lambda(\cdot,t))\,\mapsto\,
K(g_1(\cdot,t),\dots,g_\lambda(\cdot,t))
\end{equation*}
для майже всіх $t\in(0,\tau)$. Тут $g_1,\dots,g_\lambda\in L_2(\Pi)$.
Воно задає обмежений оператор
\begin{equation}\label{16f41}
K_S:\,
\bigl(H^{\sigma,\sigma/(2b);\varphi}(\Pi)\bigr)^{\lambda}
\rightarrow H^{\sigma,\sigma/(2b);\varphi}(S)
\end{equation}
для кожного $\sigma>0$ (див. \cite{Los15NK2}, доведення теореми 2).

Нехай
$$
v:=(v_0,v_1,\dots,v_{r-1})\in\bigoplus_{k=0}^{r-1}H^{s-2bk-b;\varphi}(\Gamma).
$$
Для кожного номера $k=0,\ldots,r-1$ покладемо $Lv_k:=(v_{1k},v_{2k},\dots,v_{\lambda k})$.
Розглянемо лінійне відображення
\begin{equation*}
T\,:\,v\,\mapsto\,
K_S\bigl(T_\Pi(v_{10},v_{11},\dots,v_{1,r-1}),\dots,
T_\Pi(v_{\lambda0},v_{\lambda1},\dots,v_{\lambda,r-1})
\bigr).
\end{equation*}
З обмеженості операторів
\eqref{16f39}, \eqref{16f40} і \eqref{16f41} випливає, що воно задає
обмежений оператор~\eqref{16f31}. Він є правим оберненим до оператора
\eqref{16f29}. Справді, для довільного вектора
$$
v=(v_0,v_1,\dots,v_{r-1})\in\bigoplus_{k=0}^{r-1}H^{s-2bk-b;\varphi}(\Gamma)
$$
і номера $k\in\{0,\dots,r-1\}$ маємо рівності
\begin{align*}
(R_\Gamma Tv)_k&=\biggl(R_\Gamma K_S\bigl(T_\Pi(v_{10},\dots,v_{1,r-1}),\dots,
T_\Pi(v_{\lambda0},\dots,v_{\lambda,r-1})
\bigr)\biggr)_k\\
&=K\biggl(\bigl(R_\Pi T_\Pi(v_{10},\dots,v_{1,r-1})\bigr)_k,\dots,
\bigl(R_\Pi T_\Pi(v_{\lambda0},\dots,v_{\lambda,r-1})\bigr)_k\biggr)\\
&=K\bigl((v_{10},\dots,v_{1,r-1})_k,\dots,
(v_{\lambda0},\dots,v_{\lambda,r-1})_k\bigr)=
K(v_{1k},\dots,v_{\lambda,k})\\
&=KLv_k=v_k.
\end{align*}
Тут індекс $k$ при векторі позначає $k$-ту компоненту цього вектора.
Отже, $R_\Gamma Tv=v$.

Лема 3 доведена.

Перейдемо до встановлення
інтерполяційної формули для просторів
$\mathcal{Q}^{s-2m,(s-2m)/(2b);\varphi}$.
З означення цих просторів випливає, що інтерполяційні формули для них
потрібні у випадку, коли $s\notin\{\sigma_0+r-1/2\,:\,r\in\mathbb{N}\}$.
Тому розглянемо такі інтервали:
$$
D_0=(\sigma_0,\sigma_0+1/2),\quad D_r=(\sigma_0+r-1/2,\sigma_0+r+1/2)\quad\mbox{для всіх}\quad r\in\mathbb{N}.
$$
Наступну лему доведемо в припущенні, що $s$ належить одному з цих інтервалів.

\textbf{Лема 4.} \it Нехай задано довільні $r\in\mathbb{N}$ і $\varphi\in\mathcal{M}$.
Припустимо, що $s_0,s,s_1\in D_{r-1}$ і $s_0<s<s_1$. Тоді
правильна така рівність
\begin{equation}\label{16f42}
\mathcal{Q}^{s-2m,(s-2m)/(2b);\varphi}=\,
\bigl[\mathcal{Q}^{s_0-2m,(s_0-2m)/(2b)},
\mathcal{Q}^{s_1-2m,(s_1-2m)/(2b)}\bigr]_{\psi}
\end{equation}
з точністю до еквівалентності норм.
\rm

Доведення формули \eqref{16f42}
спирається на твердження 2 (у частині інтерполяції
підпросторів) та інтерполяційну формулу для пари
соболєвських просторів $\mathcal{H}^{s_0-2m,(s_0-2m)/(2b)}$
і $\mathcal{H}^{s_1-2m,(s_1-2m)/(2b)}$.

Відмітимо, що елементи усіх просторів
$\mathcal{Q}^{\sigma-2m,(\sigma-2m)/(2b);\varphi}$, коли
$\sigma$ пробігає інтервал $D_{r-1}$, задовольняють одні і ті самі умови узгодження.
Справді, для всіх $\sigma\in D_{r-1}$ кожне число
$(\sigma-m_j-1/2-b)/(2b)$, де $j\in\{1,\dots,m\}$,
є дробовим, а його ціла частина
$$
\bigl[\frac{\sigma-m_j-1/2-b}{2b}\bigr]=
\bigl[\frac{\sigma_0+r-m_j-1-b}{2b}\bigr]=:q_{r,j},
$$
тобто не залежить від $\sigma$.
Це означає, що
у цьому випадку умовами узгодження для елементів просторів
$\mathcal{Q}^{\sigma-2m,(\sigma-2m)/(2b);\varphi}$
є рівності \eqref{16f8},
записані для всіх $k\in\{0,\dots,q_{r,j}\}$ і $j\in\{1,\dots,m\}$.

Побудуємо лінійний оператор $P^{(r)}$, який для кожного числа $\sigma\in D_{r-1}$
буде проектором соболєвського простору $\mathcal{H}^{\sigma-2m,(\sigma-2m)/(2b)}$
на його підпростір $\mathcal{Q}^{\sigma-2m,(\sigma-2m)/(2b)}$.
Нехай довільна вектор--функція
$$F:=\bigl(f,g_1,\dots,g_m,h_0,\dots,h_{\varkappa-1}\bigr)\in
\mathcal{H}^{\sigma-2m,(\sigma-2m)/(2b)}.$$
Для кожної її компоненти $g_j$, де $j\in\{1,\dots,m\}$, покладемо
\begin{equation}\label{16f43}
\begin{split}
&g^{*}_{j}:=g_j\quad\mbox{якщо}\quad q_{r,j}<0;\\
&g^{*}_{j}:=g_j+T(w_{j,0},\dots,w_{j,q_{r,j}})\quad\mbox{якщо}\quad q_{r,j}\geq0.
\end{split}
\end{equation}
Тут
\begin{equation*}
\begin{split}
&w_{j,0}=B_{j,0}(v_{0},\dots,v_{[m_j/(2b)]})\!\upharpoonright\!\Gamma-g_j\!\upharpoonright\!\Gamma,\\
&\dots\\
&w_{j,q_{r,j}}=B_{j,q_{r,j}}(v_{0},\dots,v_{[m_j/(2b)]+q_{r,j}})\!\upharpoonright\!\Gamma-
\partial_{t}^{q_{r,j}} g_j\!\upharpoonright\!\Gamma,
\end{split}
\end{equation*}
функції $v_0,v_1\dots,v_{[m_j/(2b)]+q_{r,j}}$ означені за рекурентними формулами \eqref{16f9}, а
$T$ є  оператором \eqref{16f31}, де $s:=\sigma-m_j-1/2$.
Розглянемо відображення
\begin{equation*}
P^{(r)}:\,\bigl(f,g_1,\dots,g_m,h_0,\dots,h_{\varkappa-1}\bigr)\mapsto
\bigl(f,g^*_1,\dots,g^*_m,h_0,\dots,h_{\varkappa-1}\bigr),
\end{equation*}
де $\bigl(f,g_1,\dots,g_m,h_0,\dots,h_{\varkappa-1}\bigr)
\in\mathcal{H}^{\sigma-2m,(\sigma-2m)/(2b)}$.
Оператор $P^{(r)}$ є шуканим. Справді, він є лінійним обмеженим оператором на просторі
$\mathcal{H}^{\sigma-2m,(\sigma-2m)/(2b)}$.
З його побудови випливає, що
$P^{(r)}F\in\mathcal{Q}^{\sigma-2m,(\sigma-2m)/(2b)}$ для довільного
$F\in\mathcal{H}^{\sigma-2m,(\sigma-2m)/(2b)}$.
Більш того, якщо $F\in\mathcal{Q}^{\sigma-2m,(\sigma-2m)/(2b)}$, то
$P^{(r)}F=F$. Справді, у цьому випадку правильні рівності
\eqref{16f8}. З них та з формул \eqref{16f43} випливає, що
$g_1^*=g_1,\dots,g_m^*=g_m$.

Тепер скористаємося твердженням 2.
Згідно з ним пара
$$
\bigl[\mathcal{Q}^{s_0-2m,(s_0-2m)/(2b)},
\mathcal{Q}^{s_1-2m,(s_1-2m)/(2b)}\bigr]
$$
є припустимою і правильна рівність
\begin{equation}\label{16f44}
\begin{split}
&\bigl[\mathcal{Q}^{s_0-2m,(s_0-2m)/(2b)},
\mathcal{Q}^{s_1-2m,(s_1-2m)/(2b)}\bigr]_{\psi}\\
&=\bigl[\mathcal{H}^{s_0-2m,(s_0-2m)/(2b)},
\mathcal{H}^{s_1-2m,(s_1-2m)/(2b)}\bigr]_{\psi}\cap
\mathcal{Q}^{s_0-2m,(s_0-2m)/(2b)}.
\end{split}
\end{equation}
Права частина цієї рівності є підпростором інтерполяційного простору
$$\bigl[\mathcal{H}^{s_0-2m,(s_0-2m)/(2b)},
\mathcal{H}^{s_1-2m,(s_1-2m)/(2b)}\bigr]_{\psi}.$$
На підставі твердження 3 і формул \eqref{16f19}, \eqref{16f21} і \eqref{16f22} можемо записати таке:
\begin{align*}
\bigl[&\mathcal{H}^{s_0-2m,(s_0-2m)/(2b)},
\mathcal{H}^{s_1-2m,(s_1-2m)/(2b)}\bigr]_{\psi} \\
&=\bigl[H^{s_0-2m,(s_0-2m)/(2b)}(\Omega)\oplus
\bigoplus_{j=1}^{m}H^{s_0-m_j-1/2,(s_0-m_j-1/2)/(2b)}(S)
\oplus\bigoplus_{k=0}^{\varkappa-1}H^{s_0-2bk-b}(G), \\
&\quad\;\;\, H^{s_1-2m,(s_1-2m)/(2b)}(\Omega)\oplus
\bigoplus_{j=1}^{m}H^{s_1-m_j-1/2,(s_1-m_j-1/2)/(2b)}(S)
\oplus\bigoplus_{k=0}^{\varkappa-1}H^{s_1-2bk-b}(G)\bigr]_{\psi}\\
&=\bigl[H^{s_0-2m,(s_0-2m)/(2b)}(\Omega),H^{s_1-2m,(s_1-2m)/(2b)}(\Omega)\bigr]_{\psi}\\
&\qquad
\oplus\bigoplus_{j=1}^{m}
\bigl[H^{s_0-m_j-1/2,\,(s_0-m_j-1/2)/(2b)}(S),
H^{s_1-m_j-1/2,\,(s_1-m_j-1/2)/(2b)}(S)\bigr]_{\psi}\\
&\qquad
\oplus\bigoplus_{k=0}^{\varkappa-1}\bigl[H^{s_0-2bk-b}(G),H^{s_1-2bk-b}(G)\bigr]_{\psi}\\
&=H^{s-2m,(s-2m)/(2b);\varphi}(\Omega)
\oplus\bigoplus_{j=1}^{m}H^{s-m_j-1/2,(s-m_j-1/2)/(2b);\varphi}(S)
\oplus\bigoplus_{k=0}^{\varkappa-1}H^{s-2bk-b;\varphi}(G)\\
&=\mathcal{H}^{s-2m,(s-2m)/(2b);\varphi}.
\end{align*}
Отже,
\begin{equation}\label{16f45}
\bigl[\mathcal{H}^{s_0-2m,(s_0-2m)/(2b)},
\mathcal{H}^{s_1-2m,(s_1-2m)/(2b)}\bigr]_{\psi}=
\mathcal{H}^{s-2m,(s-2m)/(2b);\varphi}
\end{equation}
з точністю до еквівалентності норм.
На підставі \eqref{16f44} і \eqref{16f45} маємо
\begin{align*}
&\bigl[\mathcal{Q}^{s_0-2m,(s_0-2m)/(2b)},
\mathcal{Q}^{s_1-2m,(s_1-2m)/(2b)}\bigr]_{\psi}\\
&=\mathcal{H}^{s-2m,(s-2m)/(2b);\varphi}\cap
\mathcal{Q}^{s_0-2m,(s_0-2m)/(2b)}=
\mathcal{Q}^{s-2m,(s-2m)/(2b);\varphi}.
\end{align*}
Остання рівність є правильною, оскільки, як було зазначено
вище, елементи просторів $\mathcal{Q}^{s_0-2m,(s_0-2m)/(2b)}$
і $\mathcal{Q}^{s-2m,(s-2m)/(2b);\varphi}$ задовольняють одні
і ті самі умови узгодження.

Лема 4 доведена.

Наостанок будуть потрібні такі дві леми, якими скористаємось при
доведенні теорем~1 і 2 у випадку $s\in\{\sigma_0+r-1/2\,:\,r\in\mathbb{N}\}$.

\textbf{Лема 5} \it
Для довільних дійсних чисел $s$ і $\varepsilon$ таких, що
$s>\varepsilon>0$, і довільного функціонального параметра
$\varphi\in\mathcal{M}$ правильна рівність
\begin{equation}\label{16f46}
H^{s,s/(2b);\varphi}(V)=
\bigl[H^{s-\varepsilon,(s-\varepsilon)/(2b);\varphi}(V),
H^{s+\varepsilon,(s+\varepsilon)/(2b);\varphi}(V)\bigr]_{1/2}
\end{equation}
з точністю до еквівалентності норм. Тут $V\in\{\Omega,S\}$.
\rm

\textbf{Доведення.}
Нехай $V=\Omega$.
Виберемо дійсне $\delta>0$ так, щоб $s-\varepsilon-\delta>0$.
Згідно з \eqref{16f22} маємо
\begin{equation*}
H^{s-\varepsilon,(s-\varepsilon)/(2b);\varphi}(\Omega)=
\bigl[H^{s-\varepsilon-\delta,(s-\varepsilon-\delta)/(2b)}(\Omega),
H^{s+\varepsilon+\delta,(s+\varepsilon+\delta)/(2b)}(\Omega)\bigr]_{\alpha}
\end{equation*}
і
\begin{equation*}
H^{s+\varepsilon,(s+\varepsilon)/(2b);\varphi}(\Omega)=
\bigl[H^{s-\varepsilon-\delta,(s-\varepsilon-\delta)/(2b)}(\Omega),
H^{s+\varepsilon+\delta,(s+\varepsilon+\delta)/(2b)}(\Omega)\bigr]_{\beta}.
\end{equation*}
Тут інтерполяційні параметри $\alpha$ і $\beta$
означаються за формулами
\begin{equation*}
\alpha(t):=t^{\delta/(2\varepsilon+2\delta)}\varphi(t^{1/(2\varepsilon+2\delta)}),
\quad
\beta(t):=t^{(2\varepsilon+\delta)/(2\varepsilon+2\delta)}\varphi(t^{1/(2\varepsilon+2\delta)})
\quad\mbox{при}\quad t\geq1
\end{equation*}
і $\alpha(t)=\beta(t):=1$ при $0<t<1$.
Звідси за твердженням 4 маємо
\begin{equation}\label{16f47}
\begin{split}
&\bigl[H^{s-\varepsilon,(s-\varepsilon)/(2b);\varphi}(\Omega),
H^{s+\varepsilon,(s+\varepsilon)/(2b);\varphi}(\Omega)\bigr]_{1/2}\\
&=\Bigl[
\bigl[H^{s-\varepsilon-\delta,(s-\varepsilon-\delta)/(2b)}(\Omega),
H^{s+\varepsilon+\delta,(s+\varepsilon+\delta)/(2b)}(\Omega)\bigr]_{\alpha},\\
&\qquad\bigl[H^{s-\varepsilon-\delta,(s-\varepsilon-\delta)/(2b)}(\Omega),
H^{s+\varepsilon+\delta,(s+\varepsilon+\delta)/(2b)}(\Omega)\bigr]_{\beta}
\Bigr]_{1/2}\\
&=\bigl[H^{s-\varepsilon-\delta,(s-\varepsilon-\delta)/(2b)}(\Omega),
H^{s+\varepsilon+\delta,(s+\varepsilon+\delta)/(2b)}(\Omega)\bigr]_{\omega}.
\end{split}
\end{equation}
Тут інтерполяційний параметр $\omega$
визначається за формулами
\begin{equation*}
\omega(t):=\alpha(t)(\beta(t)/\alpha(t))^{1/2}=t^{1/2}\varphi(t^{1/(2\varepsilon+2\delta)})
\quad\mbox{при}\quad t\geq1
\end{equation*}
і $\omega(t):=1$ при $0<t<1$.
Тому на підставі формули \eqref{16f22} маємо
\begin{equation}\label{16f48}
\bigl[H^{s-\varepsilon-\delta,(s-\varepsilon-\delta)/(2b)}(\Omega),
H^{s+\varepsilon+\delta,(s+\varepsilon+\delta)/(2b)}(\Omega)\bigr]_{\omega}=
H^{s,s/(2b);\varphi}(\Omega).
\end{equation}
Тепер з рівностей \eqref{16f47} і \eqref{16f48} випливає формула \eqref{16f46}.
Оскільки рівності \eqref{16f47} і \eqref{16f48} виконуються з точністю
до еквівалентності норм, то так само буде і для рівності \eqref{16f46}.

Для випадку $V=S$ доведення проводиться точно так, тільки замість \eqref{16f22}
використовуємо формулу \eqref{16f21}.

Лема 5 доведена.

\textbf{Лема 6} \it
Для довільних додатних чисел $s$ і $\varepsilon$ таких, що
$s-\varepsilon>\sigma_0$, і довільного функціонального параметра
$\varphi\in\mathcal{M}$ правильна рівність
\begin{equation}\label{16f67}
\mathcal{H}^{s-2m,(s-2m)/(2b);\varphi}=
\bigl[\mathcal{H}^{s-2m-\varepsilon,(s-2m-\varepsilon)/(2b);\varphi},
\mathcal{H}^{s-2m+\varepsilon,(s-2m+\varepsilon)/(2b);\varphi}\bigr]_{1/2}
\end{equation}
з точністю до еквівалентності норм.
\rm

\textbf{Доведення.}
Рівність \eqref{16f67} випливає з твердження~3, формул \eqref{16f46} та
їх аналогу для ізотропних просторів $H^{s;\varphi}(G)$ (див.\cite[Лема~4.3]{MikhailetsMurach14}).
Справді,
\begin{align*}
\bigl[&\mathcal{H}^{s-2m-\varepsilon,(s-2m-\varepsilon)/(2b);\varphi},
\mathcal{H}^{s-2m+\varepsilon,(s-2m+\varepsilon)/(2b);\varphi}\bigr]_{1/2} \\
&=\bigl[H^{s-2m-\varepsilon,(s-2m-\varepsilon)/(2b);\varphi}(\Omega)\oplus
\bigoplus_{j=1}^{m}H^{s-m_j-1/2-\varepsilon,(s-m_j-1/2-\varepsilon)/(2b);\varphi}(S)\\
&\qquad\oplus\bigoplus_{k=0}^{\varkappa-1}H^{s-2bk-b-\varepsilon;\varphi}(G), \\
&\quad\;\;\, H^{s-2m+\varepsilon,(s-2m+\varepsilon)/(2b);\varphi}(\Omega)\oplus
\bigoplus_{j=1}^{m}H^{s-m_j-1/2+\varepsilon,(s-m_j-1/2+\varepsilon)/(2b);\varphi}(S)\\
&\qquad\oplus\bigoplus_{k=0}^{\varkappa-1}H^{s-2bk-b+\varepsilon;\varphi}(G)\bigr]_{1/2}\\
&=\bigl[H^{s-2m-\varepsilon,(s-2m-\varepsilon)/(2b);\varphi}(\Omega),
H^{s-2m+\varepsilon,(s-2m+\varepsilon)/(2b);\varphi}(\Omega)\bigr]_{1/2}\\
&\qquad
\oplus\bigoplus_{j=1}^{m}
\bigl[H^{s-m_j-1/2-\varepsilon,\,(s-m_j-1/2-\varepsilon)/(2b);\varphi}(S),
H^{s-m_j-1/2+\varepsilon,\,(s-m_j-1/2+\varepsilon)/(2b);\varphi}(S)\bigr]_{1/2}\\
&\qquad
\oplus\bigoplus_{k=0}^{\varkappa-1}\bigl[H^{s-2bk-b-\varepsilon;\varphi}(G),
H^{s-2bk-b+\varepsilon;\varphi}(G)\bigr]_{1/2}\\
&=H^{s-2m,(s-2m)/(2b);\varphi}(\Omega)
\oplus\bigoplus_{j=1}^{m}H^{s-m_j-1/2,(s-m_j-1/2)/(2b);\varphi}(S)
\oplus\bigoplus_{k=0}^{\varkappa-1}H^{s-2bk-b;\varphi}(G)\\
&=\mathcal{H}^{s-2m,(s-2m)/(2b);\varphi}.
\end{align*}
Отже, правильно \eqref{16f67}.

Лема 6 доведена.

\textbf{6. Доведення основних результатів.}

У цьому пункті доведемо теореми 1 і 2.

\textbf{Доведення теореми 1.}
Спочатку розглянемо випадок $s\notin\{\sigma_0+r-1/2\,:\,r\in\mathbb{N}\}$.
Тоді $s\in D_{r-1}$ для деякого $r\in\mathbb{N}$.
Виберемо такі числа $s_0,s_1\in D_{r-1}$,
що $s_0<s<s_1$.
Завдяки згаданій вище теоремі М.~В.~Житарашу \cite[теорема 9.1]{Zhitarashu85}
маємо ізоморфізми у просторах Соболева
\begin{equation}\label{16f49}
\Lambda:\,H^{s_j,s_j/(2b)}(\Omega)\leftrightarrow
\mathcal{Q}^{s_j-2m,(s_j-2m)/(2b)}
\quad\mbox{для кожного}\quad j\in\{0,1\}.
\end{equation}
Застосувавши інтерполяцію з функціональним параметром $\psi$
до \eqref{16f49}, отримаємо ще один ізоморфізм
\begin{equation}\label{16f50}
\Lambda:\,
\bigl[H^{s_0,s_0/(2b)}(\Omega),
H^{s_1,s_1/(2b)}(\Omega)\bigr]_{\psi}\leftrightarrow
[\mathcal{Q}^{s_0-2m,(s_0-2m)/(2b)},
\mathcal{Q}^{s_1-2m,(s_1-2m)/(2b)}]_{\psi}.
\end{equation}
Цей ізоморфізм є розширенням по неперервності
відображення \eqref{16f4} оскільки $C^{\infty}(\overline{\Omega})$
щільно в області визначення \eqref{16f50}.
Застосувавши у \eqref{16f50} інтерполяційні формули
\eqref{16f22} та \eqref{16f42} отримаємо \eqref{16f11}.

Нехай тепер $s\in\{\sigma_0+r-1/2\,:\,r\in\mathbb{N}\}$.
За доведеним маємо ізоморфізми
\begin{equation}\label{16f51}
\Lambda:\,H^{s\pm\varepsilon,(s\pm\varepsilon)/(2b);\varphi}(\Omega)\leftrightarrow
\mathcal{Q}^{s-2m\pm\varepsilon,(s-2m\pm\varepsilon)/(2b);\varphi}.
\end{equation}
Тут $\varepsilon\in(0;1/2)$ те саме, що і у \eqref{16f10}.
Застосувавши інтерполяцію з числовим параметром $1/2$
до \eqref{16f51}, отримаємо ще один ізоморфізм
\begin{equation}\label{16f52}
\begin{split}
\Lambda:\,
&\bigl[H^{s-\varepsilon,(s-\varepsilon)/(2b);\varphi}(\Omega),
H^{s+\varepsilon,(s+\varepsilon)/(2b);\varphi}(\Omega)\bigr]_{1/2}\leftrightarrow\\
&\bigl[\mathcal{Q}^{s-2m-\varepsilon,(s-2m-\varepsilon)/(2b);\varphi},
 \mathcal{Q}^{s-2m+\varepsilon,(s-2m+\varepsilon)/(2b);\varphi}\bigr]_{1/2}.
\end{split}
\end{equation}
Тепер застосуємо у \eqref{16f52} інтерполяційну формулу \eqref{16f46} і
означення \eqref{16f10} простору $\mathcal{Q}^{s-2m,(s-2m)/(2b);\varphi}$.
В результаті отримаємо ізоморфізм \eqref{16f11}.

Теорема 1 доведена.

\textbf{Зауваження 1.} \it Означений за формулою \eqref{16f10} простір
не залежать з точністю до еквівалентності норм від
вибору числа $\varepsilon\in(0,1/2)$.
\rm

Справді, з теореми~1 випливає, що для $s\in\{\sigma_0+r-1/2\,:\,r\in\mathbb{N}\}$
виконується ізоморфізм
\begin{equation}\label{16f53}
\Lambda:\,H^{s,s/(2b);\varphi}(\Omega)\leftrightarrow
\mathcal{Q}_{\,\varepsilon}^{s-2m,(s-2m)/(2b);\varphi}.
\end{equation}
Тут
\begin{equation*}
\mathcal{Q}_{\,\varepsilon}^{s-2m,(s-2m)/(2b);\varphi}:=\,
\bigl[\mathcal{Q}^{s-2m-\varepsilon,(s-2m-\varepsilon)/(2b);\varphi},
\mathcal{Q}^{s-2m+\varepsilon,(s-2m+\varepsilon)/(2b);\varphi}\bigr]_{1/2}.
\end{equation*}
З ізоморфізму \eqref{16f53} негайно випливає, що простір
$\mathcal{Q}_{\,\varepsilon}^{s-2m,(s-2m)/(2b);\varphi}$
не залежить від $\varepsilon$ з точністю до еквівалентності норм.

\textbf{Доведення теореми 2.}
Спочатку покажемо, що з умов \eqref{16f13}--\eqref{16f15} цієї теореми випливає правильність імплікації
\begin{equation}\label{16f54}
u\in H^{\sigma-\lambda,(\sigma-\lambda)/(2b);\varphi}_{\mathrm{loc}}
(\omega,\pi_1)\;\Rightarrow\;u\in H^{\sigma-\lambda+1,(\sigma-\lambda+1)/(2b);\varphi}_{\mathrm{loc}}
(\omega,\pi_1)
\end{equation}
для кожного цілого $\lambda\geq1$, такого що $\sigma-\lambda+1>\sigma_{0}$.

Виберемо довільну функцію $\chi\in C^\infty(\overline\Omega)$ з $\mbox{supp}\,\chi\subset\omega\cup\pi_1$.
Для $\chi$ існує функція $\eta\in C^\infty(\overline\Omega)$ така, що $\mbox{supp}\,\eta\subset\omega\cup\pi_1$ і
$\eta=1$ в околі $\mbox{supp}\,\chi$. Переставляючи кожний з диференцiальних операторiв
$A$,  $B_{j}$ і $\partial^k_t$ з оператором множення на  $\chi$, можемо записати
\begin{equation}\label{16f55}
\begin{split}
\Lambda(\chi u)&=\Lambda(\chi\eta u)=\chi\,\Lambda(\eta u)+ \Lambda'(\eta u)\\
&=\chi\,\Lambda u+\Lambda'(\eta u)=\chi\,(f,g_{1},\dots,g_{m},h_{0},\dots,h_{\varkappa-1})+\Lambda'(\eta u).
\end{split}
\end{equation}
Тут $\Lambda':=(A',B'_{1},\ldots,B'_{m},C'_{0},\dots C'_{\varkappa-1})$-- диференціальний оператор з компонентами
\begin{equation}\label{16f56}
A'(x,t,D_x,\partial_t)=\sum_{|\alpha|+2b\beta\leq 2m-1}a^{\alpha,\beta}_{1}(x,t)\,D^\alpha_x\partial^\beta_t,
\end{equation}
\begin{equation}\label{16f57}
B_{j}'(x,t,D_x,\partial_t)=\sum_{|\alpha|+2b\beta\leq m_j-1}
b_{j,1}^{\alpha,\beta}(x,t)\,D^\alpha_x\partial^\beta_t,\quad j=1,\ldots,m
\end{equation}
і
\begin{equation}\label{16f58}
C_{0}'=0,\quad
C_{k}'(x,t,\partial_t)=\sum_{l=0}^{k-1}
c_{l,\,k}(x,t)\,\partial^{\,l}_t,\quad k=1,\ldots,\varkappa-1,
\end{equation}
де всі $a^{\alpha,\beta}_{1}\in C^{\infty}(\overline{\Omega})$,
$b_{j,1}^{\alpha,\beta}\in C^{\infty}(\overline{S})$ і $c_{l,\,k}\in C^{\infty}(\overline{G})$.
Цей оператор неперервно діє у парі просторів
\begin{equation}\label{16f59}
\Lambda':\,H^{s,s/(2b);\varphi}(\Omega)\rightarrow
\mathcal{H}^{s+1-2m,(s+1-2m)/(2b);\varphi}
\end{equation}
для кожного $s>\sigma_{0}-1$. У випадку, коли $\varphi\equiv1$, це зразу слідує з  \eqref{16f56}, \eqref{16f57},
\eqref{16f58} і відомих властивостей анізотропного простору Соболєва $H^{s,s/(2b)}(\Omega)$ (див, наприклад, \cite[гл.~I, лема~4, та гл.~II, теореми~3 і~7]{Slobodetskii58}).
Обмеженість оператора \eqref{16f59} у загальній ситуації зразу випливає з цього випадку за допомогою
інтерполяційних формул \eqref{16f19}, \eqref{16f21}, \eqref{16f22}.

З умов \eqref{16f13}, \eqref{16f14} та \eqref{16f15} маємо включення
$$
\chi\,(f,g_{1},\dots,g_{m},h_0,\dots,h_{\varkappa-1})
\in\mathcal{H}^{\sigma-2m,(\sigma-2m)/(2b);\varphi}.
$$
Крім того, згідно \eqref{16f59} з $s:=\sigma-\lambda$, має місце імплікація
\begin{equation*}
\begin{split}
u&\in H^{\sigma-\lambda,(\sigma-\lambda)/(2b);\varphi}_{\mathrm{loc}}
(\omega,\pi_{1})\\
&\Rightarrow\;\Lambda'(\eta u)\in
\mathcal{H}^{\sigma-\lambda+1-2m,(\sigma-\lambda+1-2m)/(2b);\varphi}.
\end{split}
\end{equation*}
Тому, скориставшись \eqref{16f55} можемо записати
\begin{equation}\label{16f60}
\begin{split}
u&\in H^{\sigma-\lambda,(\sigma-\lambda)/(2b);\varphi}_{\mathrm{loc}}
(\omega,\pi_{1})\\
&\Rightarrow\;\Lambda(\chi u)\in
\mathcal{H}^{\sigma-\lambda+1-2m,(\sigma-\lambda+1-2m)/(2b);\varphi}.
\end{split}
\end{equation}
Тепер покажемо, що для довільного $s>\sigma_0$ з включення
$\Lambda(\chi u)\in
\mathcal{H}^{s-2m,(s-2m)/(2b);\varphi}$ і зробленого вибору функції $\chi$ випливає
включення
\begin{equation}\label{16f61}
\Lambda(\chi u)\in
\mathcal{Q}^{s-2m,(s-2m)/(2b);\varphi}.
\end{equation}
Попередньо відмітимо, що
оскільки $\mathrm{dist}(\mbox{supp}\,\chi,\Gamma)>0$, то $\chi=0$ в деякому околі $\Gamma$.
Звідси маємо, що
\begin{equation}\label{16f62}
\Lambda(\chi u)=0
\end{equation}
в цьому околі $\Gamma$. Розглянемо окремо випадки
$s\notin\{\sigma_0+r-1/2\,:\,r\in\mathbb{N}\}$ і $s\in\{\sigma_0+r-1/2\,:\,r\in\mathbb{N}\}$.

Нехай $s\notin\{\sigma_0+r-1/2\,:\,r\in\mathbb{N}\}$. Рівність \eqref{16f62}
означає, що для
вектор--функції $\Lambda(\chi u)\in
\mathcal{H}^{s-2m,(s-2m)/(2b);\varphi}$
виконуються умови узгодження \eqref{16f8}, і отже, правильне включення
\eqref{16f61}.

Нехай $s\in\{\sigma_0+r-1/2\,:\,r\in\mathbb{N}\}$.
Виберемо функцію $\chi_1\in C^{\infty}(\overline\Omega)$ таку, що $\chi_1=0$ в деякому
околі $\Gamma$. Зі сказаного у попередньому абзаці випливає, що відображення
$M_{\chi_1}\,:\,F\mapsto\chi_1 F$ визначає обмежені оператори
\begin{equation}\label{16f63}
M_{\chi_1}\,:\,\mathcal{H}^{s-2m\pm\varepsilon,(s-2m\pm\varepsilon)/(2b);\varphi}
\mapsto\mathcal{Q}^{s-2m\pm\varepsilon,(s-2m\pm\varepsilon)/(2b);\varphi}
\end{equation}
для довільного $\varepsilon\in(0,1/2)$.
Застосуємо до \eqref{16f63} інтерполяцію з числовим параметром $1/2$ і
отримаємо ще один обмежений оператор
\begin{equation}\label{16f64}
\begin{split}
M_{\chi_1}\,:\,
&\bigl[\mathcal{H}^{s-2m-\varepsilon,(s-2m-\varepsilon)/(2b);\varphi},
\mathcal{H}^{s-2m+\varepsilon,(s-2m+\varepsilon)/(2b);\varphi}\bigr]_{1/2}
\to \\
&\bigl[\mathcal{Q}^{s-2m-\varepsilon,(s-2m-\varepsilon)/(2b);\varphi},
\mathcal{Q}^{s-2m+\varepsilon,(s-2m+\varepsilon)/(2b);\varphi}\bigr]_{1/2}.
\end{split}
\end{equation}
Тепер застосуємо у \eqref{16f64} інтерполяційну формулу  \eqref{16f67} і
означення \eqref{16f10} простору $\mathcal{Q}^{s-2m,(s-2m)/(2b);\varphi}$.
В результаті отримаємо обмежений оператор
\begin{equation}\label{16f65}
M_{\chi_1}\,:\,
\mathcal{H}^{s-2m,(s-2m)/(2b);\varphi}\mapsto
\mathcal{Q}^{s-2m,(s-2m)/(2b);\varphi}.
\end{equation}

Нехай $\Lambda(\chi u)\in\mathcal{H}^{s-2m,(s-2m)/(2b);\varphi}$. Позначимо
$O_1$ той окіл $\Gamma$, у якому $\Lambda(\chi u)=0$.
Нехай $O_2$ і $O_3$ такі околи $\Gamma$, що $O_3\subset O_2\subset O_1$,
і $\overline{O_3}\subset O_2$, і $\overline{O_2}\subset O_1$.
Виберемо функцію $\chi_1\in C^{\infty}(\overline\Omega)$ таку, що
$\chi_1=0$ в $O_3$ і $\chi_1=1$ в $\overline\Omega\setminus O_2$.
Тоді $\chi_1\Lambda(\chi u)=\Lambda(\chi u)$. А з \eqref{16f65}
випливає, що $\chi_1\Lambda(\chi u)\in\mathcal{Q}^{s-2m,(s-2m)/(2b);\varphi}$.
Отже, правильне включення \eqref{16f61}.

З імплікації \eqref{16f60}, включення \eqref{16f61} та наслідку~1 випливає, що
\begin{equation*}
\begin{split}
u&\in H^{\sigma-\lambda,(\sigma-\lambda)/(2b);\varphi}_{\mathrm{loc}}
(\omega,\pi_{1})\\
&\Rightarrow\;\Lambda(\chi u)\in
\mathcal{Q}^{\sigma-\lambda+1-2m,(\sigma-\lambda+1-2m)/(2b);\varphi}
\\
&\Rightarrow\;\chi u\in
H^{\sigma-\lambda+1,(\sigma-\lambda+1)/(2b);\varphi}(\Omega).
\end{split}
\end{equation*}
Відмітимо, що тут наслідок~1 застосовний, оскільки $\chi u\in H^{\sigma_0,\sigma_0/(2b)}(\Omega)$
за умовою теореми, і $\sigma-\lambda+1>\sigma_{0}$.
Тим самим, імплікація \eqref{16f54} доведена, якщо зважити на зроблений вибір функції $\chi$.

Використаємо цю імплікацію для доведення включення
$u\in H^{\sigma,\sigma/(2b);\varphi}_{\mathrm{loc}}(\omega,\pi_1)$.
Розглянемо окремо випадки $\sigma\notin\mathbb{Z}$ і $\sigma\in\mathbb{Z}$.

Нехай спочатку $\sigma\notin\mathbb{Z}$. У цьому випадку існує ціле число $\lambda_{0}\geq1$ таке, що
\begin{equation}\label{16f66}
\sigma-\lambda_{0}<\sigma_{0}<\sigma-\lambda_{0}+1.
\end{equation}
Скориставшись імплікацією \eqref{16f54}
послідовно для значень $\lambda:=\lambda_{0}$,
$\lambda:=\lambda_{0}-1$,..., $\lambda:=1$, виводимо необхідне включення слідуючим чином:
\begin{equation*}
\begin{split}
u&\in H^{\sigma_0,\sigma_0/(2b)}(\Omega)\subset
H^{\sigma-\lambda_{0},(\sigma-\lambda_{0})/(2b);\varphi}_{\mathrm{loc}}(\omega,\pi_1)\\
&\Rightarrow\;u\in
H^{\sigma-\lambda_{0}+1,(\sigma-\lambda_{0}+1)/(2b);\varphi}_{\mathrm{loc}}(\omega,\pi_1)\\
&\Rightarrow\;\ldots\;\Rightarrow\;u\in H^{\sigma,\sigma/(2b);\varphi}_{\mathrm{loc}}(\omega,\pi_1).
\end{split}
\end{equation*}
Відмітимо, що $u\in H^{\sigma_0,\sigma_0/(2b)}(\Omega)$ за умовою теореми.

Нехай тепер $\sigma\in\mathbb{Z}$. У цьому випадку не існує цілого числа $\lambda_{0}$,
що задовольняє~\eqref{16f66}.
Але, оскільки $\sigma-1/4\notin\mathbb{Z}$ і $\sigma-1/4>\sigma_{0}$, то, як довели у попередньому
абзаці, правильне включення
\begin{equation*}
u\in H^{\sigma-1/4,(\sigma-1/4)/(2b);\varphi}_{\mathrm{loc}}(\omega,\pi_1).
\end{equation*}
Звідси, скориставшись імплікацією \eqref{16f54} з  $\lambda:=1$, виводимо потрібне включення, а саме:
\begin{align*}
u\in H^{\sigma-1/4,(\sigma-1/4)/(2b);\varphi}_{\mathrm{loc}}(\omega,\pi_1)\subset
H^{\sigma-1,(\sigma-1)/(2b);\varphi}_{\mathrm{loc}}(\omega,\pi_1)
\;\Rightarrow\;u\in H^{\sigma,\sigma/(2b);\varphi}_{\mathrm{loc}}(\omega,\pi_1).
\end{align*}

\end{document}